\documentclass[11pt]{amsart}
\usepackage{amsxtra}
\usepackage{amssymb}
\theoremstyle{plain}

\newcommand{\cleqn}{\setcounter{equation}{0}}
\newcommand{\clth}{\setcounter{theorem}{0}}
\newcommand {\sectionnew}[1]{\section{#1}\cleqn\clth}

\newtheorem{theorem}{Theorem}[section]
\newtheorem{lemma}[theorem]{Lemma}
\newtheorem{definition-theorem}[theorem]{Definition-Theorem}
\newtheorem{proposition}[theorem]{Proposition}
\newtheorem{corollary}[theorem]{Corollary}
\newtheorem{definition}[theorem]{Definition}
\newtheorem{example}[theorem]{Example}

\newtheorem{remark}[theorem]{Remark}
\newtheorem{conjecture}[theorem]{Conjecture}
\newtheorem{notation}[theorem]{Notation}

\newcommand \bth[1] { \begin{theorem}\label{t#1} }
\newcommand \ble[1] { \begin{lemma}\label{l#1} }

\newcommand \bpr[1] { \begin{proposition}\label{p#1} }
\newcommand \bco[1] { \begin{corollary}\label{c#1} }
\newcommand \bde[1] { \begin{definition}\label{d#1}\rm }
\newcommand \bex[1] { \begin{example}\label{e#1}\rm }
\newcommand \bre[1] { \begin{remark}\label{r#1}\rm }
\newcommand \bcj[1] { \begin{conjecture}\label{j#1}\rm }

\newcommand \bnota[1] { \begin{notation}\label{n#1}\rm }
\renewcommand {\eth} { \end{theorem} }
\newcommand {\ele} { \end{lemma} }

\newcommand {\epr} { \end{proposition} }
\newcommand {\eco} { \end{corollary} }
\newcommand {\ede} { \end{definition} }
\newcommand {\eex} { \end{example} }
\newcommand {\ere} { \end{remark} }
\newcommand {\ecj} { \end{conjecture} }

\newcommand {\enota} { \end{notation} }
\newcommand \thref[1]{Theorem \ref{t#1}}
\newcommand \leref[1]{Lemma \ref{l#1}}
\newcommand \prref[1]{Proposition \ref{p#1}}
\newcommand \coref[1]{Corollary \ref{c#1}}

\newcommand \deref[1]{Definition \ref{d#1}}

\newcommand \reref[1]{Remark \ref{r#1}}

\def \Cset {{\mathbb C}}

\def \Zset {{\mathbb Z}}

\def \Nset {{\mathbb N}}

\def \QP {{\mathcal{QP}}}
\def \QA {{\mathcal{QA}}}

\def \OO {{\mathcal{O}}}

\def \al {\alpha}

\def \la {\lambda}
\def \La {\Lambda}

\def \Om {\Omega}

\def \Sig {\Sigma}
\def \sig {\sigma}

\def \vph {\varphi}

\def \sig{\sigma}

\def \mt  {\mapsto}

\def \ol {\overline}
\def \wt {\widetilde}

\DeclareMathOperator \vker { \operatorname{vker} }

\DeclareMathOperator \Ann { {\mathrm{Ann}} }

\DeclareMathOperator \Ai {\operatorname{Ai}}

\begin{document}
\title[Noncommutative Bispectral Darboux Transformations]
{Noncommutative Bispectral Darboux Transformations}
\author[Joel Geiger]{Joel Geiger}
\address{
Department of Mathematics\\
Massachusetts Institute of Technology\\
Cambridge, MA 02139-4307 \\
U.S.A.
}
\email{jbgeiger@math.mit.edu}
\author[Emil Horozov]{Emil Horozov}
\address{
Department of Mathematics and Informatics \\
Sofia University \\
5 J. Bourchier Blvd. \\
Sofia 1126 \\
and 
Institute of Mathematics and Informatics \\ 
                          Bulg. Acad. of Sci., Acad. G. Bonchev Str., Block 8, 1113 Sofia \\
Bulgaria
}
\email{horozov@fmi.uni-sofia.bg}
\author[Milen Yakimov]{Milen Yakimov}
\thanks{The research of M.Y. was partially supported by NSF grant DMS-1303036
and Louisiana Board of Regents grant Pfund-403.}
\address{
Department of Mathematics \\
Louisiana State University \\
Baton Rouge, LA 70803 \\
U.S.A.
}
\email{yakimov@math.lsu.edu}
\date{}
\keywords{The noncommutative bispectral problem, Darboux transformations, matrix rank one bispectral functions, 
the Airy bispectral function}
\subjclass[2010]{Primary 37K35; Secondary 16S32, 39A70}
\begin{abstract} 
We prove a general theorem establishing the bispectrality of noncommutative 
Darboux transformations. It has a wide range of applications that establish bispectrality 
of such transformations for differential, difference and $q$-difference operators with 
values in all noncommutative algebras. All known bispectral Darboux transformations
are special cases of the theorem. Using the methods of quasideterminants and the spectral 
theory of matrix polynomials, we explicitly classify the set of bispectral Darboux transformations 
from rank one differential operators and Airy operators with values in matrix algebras. These sets 
generalize the classical Calogero--Moser spaces and Wilson's adelic Grassmannian.
\end{abstract}
\maketitle
\sectionnew{Introduction}
\label{intro}
\subsection{}
\label{1.1}
The study of bispectral operators was initiated by Duistermaat and Gr\"unbaum in \cite{DG} motivated by their applications to computer tomography and time-band limiting. Since 1985 the bispectral problem has attracted a lot of attention in pure mathematics and has been related to many diverse areas, such as integrable systems (KP hierarchy, Sato's Grassmannian and the Calogero-Moser system) \cite{BHY2,DG,KR,Wi1,Wi2}, orthogonal polynomials \cite{GH1,GH2}, 
representation theory ($W_{1+\infty}$, Virasoro and Kac--Moody algebras) \cite{BHY3,FMTV}, 
ideal structure and automophisms of the first Weyl algebra 
\cite{BHY1,BW} and many others. Let $\Om_1$ and $\Om_2$ be two domains in $\Cset$.
The continuous-continuous (scalar) bispectral problem asks for finding all analytic functions
\[
\Psi \colon \Om_1 \times \Om_2 \to \Cset
\]
for which there exist analytic differential operators 
$L(x, \partial_x)$ and $\La(z, \partial_z)$ on 
$\Om_1$ and $\Om_2$, and two analytic functions $\theta \colon \Om_1 \to \Cset$, $f \colon \Om_2 \to \Cset$,
such that
\begin{align}
L(x, \partial_x) \Psi(x,z) &= f(z) \Psi(x,z),  
\label{eq1}
\\
\theta(x) \Psi(x, z) &=
\La(z, \partial_z) \Psi(x,z) 
\label{eq2}
\end{align}
on $\Om_1 \times \Om_2$. There are discrete, $q$-difference 
and mixed versions of the problem \cite{GH1,GH2, HI}, which have played a similarly important role.

The noncommutative version of the problem is stated for a (generally noncommutative) associative complex finite dimensional algebra 
$R$. It asks for classifying all $R$-valued analytic functions
\[
\Psi \colon \Om_1 \times \Om_2 \to R
\]
for which there exist $R$-valued analytic differential operators 
$L(x, \partial_x)$ and $\La(z, \partial_z)$ on 
$\Om_1$ and $\Om_2$, and $R$-valued analytic functions $\theta \colon \Om_1 \to R$, 
$f \colon \Om_2 \to R$, such that
\begin{align}
L(x, \partial_x) \Psi(x,z) &= f(z) \Psi(x,z), 
\label{eq3}
\\
\theta(x) \Psi(x, z) &=
\Psi(x,z) \La(z, \partial_z) 
\label{eq4} 
\end{align}
on $\Om_1 \times \Om_2$. The right action of differential operators in the second equality is defined by
\[
\Psi(x,z) \cdot  ( a(z) \partial_z^k) := (-1)^k \partial_z^k \left( \Psi(x,z) a(z) \right)
\]
and linearity. It was observed in \cite{Du,GI,LM} that the right action was needed in order to obtain nontrivial 
examples which were absent if one insisted on left actions only. More precisely, if both actions are left, then they do not commute 
with each other because the algebra $R$ is noncommutative in general. The left action in $x$ and the right action in $z$ commute with each other, thus 
giving rise to a  bimodule for two noncommutative algebras of differential operators.

The matrix bispectral problem is the problem in which one takes $R = M_n(\Cset)$ in the above 
definition. It was first considered in \cite{Z1,Z2} with relation to the ZS-AKNS operators.
In the last 20 years there has been a great deal of research 
on the continuous and discrete versions of the matrix bispectral problem and numerous relations 
were found with matrix orthogonal polynomials \cite{DuG,GPT2, M}, spherical functions \cite{GPT}, 
and many other topics. 

The construction and classification of noncommutative (or even just matrix) bispectral functions is 
significantly harder than the scalar problem. One of the major methods in the scalar case is the 
use of Darboux transformations for obtaining large families of complicated examples from simpler 
ones. In this approach, the first equation \eqref{eq1} for the transformed function 
is direct while the second one is highly nontrivial. 
This was carried out for second order operators in \cite{DG}, all rank 1 bispectral operators \cite{Wi1} in Wilson's  
classification, Bessel functions \cite{BHY2}, Airy functions \cite{BHY2,KR}, orthogonal polynomials 
\cite{GH1}
and many other situations. Even though these results are of similar nature, their proofs relied on  a variety of different algebraic, algebro-geometric and integrable systems 
techniques. A general theorem for bispectrality of scalar Darboux transformations was
obtained in \cite{BHY1}. In some cases the algebras of dual bispectral operators were also classified 
\cite{Ha,I,Wi2}. 
\subsection{}
\label{1.2}
In this paper we obtain two major sets of results. The first one is a very general theorem proving that 
Darboux transformations from noncommutative bispectral functions produce again noncommutative bispectral 
functions. All previous examples that we are aware of can be recovered as special cases 
of this theorem, for example all functions in \cite{BL,G1,Ka,Z1,Z2}. All results on bispectrality of scalar Darboux transformations \cite{BHY1,BHY2,DG,GH1,HI,I,KR,Wi1} also arise as special cases of this theorem.

Our theorem establishes a general bispectrality type result for all Darboux transformations in bimodules of 
noncommutative algebras $B$. These noncommutative algebras $B$ can be taken to be the algebra of analytic differential 
operators with values in a finite dimensional algebra $R$, the algebras of difference or $q$-difference operators with 
values in $R$, and various other examples of a similar nature. For the sake of simplicity in the introduction we 
formulate only the special case of the theorem when the algebra $B$ is the first one in the list. We refer the 
reader to \thref{GenTh} for the general result for bispectral Darboux transformations in bimodules of noncommutative 
algebras.

Let $\Psi \colon \Om_1 \times \Om_2 \to R$ be a noncommutative bispectral function as in \eqref{eq3}--\eqref{eq4}. 
We will assume 
that it is a nonsplit function in the two variables, in the sense that there are no nonzero analytic differential 
operators $L(x, \partial_x)$ and $\La(z, \partial_z)$ that satisfy \eqref{eq3} or \eqref{eq4} 
with $f(z) \equiv 0$ or $\theta(x) \equiv 0$.

To formulate our result, we need to introduce several algebras associated to the function
$\Psi(x,z)$.
Let $B_1$ be the algebra of all $R$-valued analytic differential operators $P(x, \partial_x)$ 
on $\Om_1$ for which there exits an $R$-valued analytic differential operator $S(z, \partial_z)$ on $\Om_2$ such that 
\[
P(x, \partial_x) \Psi(x,z) = \Psi(x,z) S(z, \partial_z).
\]
Let $B_2$ be the algebra consisting of all such operators $S(z, \partial_z)$. The nonsplit condition implies 
that the map
\[
b \colon B_1 \to B_2, \quad \mbox{given by} \quad b(P(x, \partial_x)):= S(z, \partial_z),
\]
is a well defined algebra isomorphism. Denote by $M$ the $B_1 - B_2$ bimodule of 
$R$-valued analytic functions on $\Om_1 \times \Om_2$.
Let $K_i$ be the subalgebras of $B_i$ consisting of $R$-valued analytic functions, 
i.e. the differential operators in $B_i$ of order $0$. Let
\[
A_1 := b^{-1}(K_2), \quad A_2 := b(K_1).
\]
These are the algebras of bispectral operators corresponding to $\Psi(x,z)$. They consist 
of the differential operators that have the properties 
\eqref{eq3}--\eqref{eq4} for some functions $f(z)$, $\theta(x)$.  
In particular,  
\[
L(x, \partial_x) \in A_1 \quad \mbox{and} \quad 
\La(z, \partial_z) \in A_2.
\]

\bth{th-i} In the above setting, for every factorization, 
\[
L(x, \partial_x) = Q(x, \partial_x) g(x)^{-1} P(x, \partial_x)
\]
with $L(x, \partial_x) \in A_1$, $P(x, \partial_x), Q(x, \partial_x) \in B_1$, 
and a regular element $g(x) \in K_1$, 
the $R$-valued function 
\[
\Phi(x,z) := P(x,\partial_x) \Psi(x,z) \colon \Om_1 \times \Om_2 \to R
\]
is bispectral, and, more precisely, satisfies
\begin{align*}
\left( P(x, \partial_x) Q(x, \partial_x)  g(x)^{-1} \right) \Phi(x,z) 
&= \Phi(x,z) b(L)(z) \\
g(x) \Phi(x,z) &= \Phi(x,z) \left( b(L)(z)^{-1} b(Q)(z, \partial_z) b(P)(z, \partial_z) \right).
\end{align*} 
\eth
We recall that a regular element of an algebra is one that is not left or right zero divisor. 

The general version of the above result (given in \thref{GenTh}) deals with arbitrary noncommutative algebras 
$B_1, B_2$ and $B_1 - B_2$ bimodules $M$. It presents a general approach 
to establish bispectrality of Darboux transformations 
in bimodules of noncommutative algebras. 

Previously, combinations of algebraic, algebro-geometric and integrable systems 
techniques were used to obtain bispectrality on a case by case basis. This theorem, not only proves all these results with a single method, but also has a wide range of new applications. 
The theorem directly relates the method of Darboux transformations from integrable systems to noncommutative algebra. We expect that it will be also helpful in relating bispectrality to matrix factorization \cite{E}.
\subsection{}
\label{1.3}
The second set of results in the paper classifies the bispectral Darboux transformations in 
\thref{th-i} in several important cases. The rank 1 bispectral functions of Wilson \cite{Wi1} have played a key role
in the classification of scalar bispectral functions. They were first introduced in \cite{Wi1} in terms of wave functions 
for the KP hierarchy with unicursal spectral curves. In \cite{BHY2} it was proved that they are precisely the set of all
bispectral Darboux transformations from the function $\Psi(x,z) = e^{xz}$. Our first result 
classifies all bispectral Darboux transformations from the function
\[
\Psi(x,z) = e^{xz} I_n \quad \mbox{in the matrix case} \quad R = M_n(\Cset).
\]
In this case $\Om_1 = \Om_2 = \Cset$, and $B_1$ and $B_2$ are the algebras of matrix differential operators 
with polynomial coefficients in the variables $x$ and $z$, respectively. The isomorphism $b$ is given by 
\[
b(\partial_x) = z, \; b(x) = - \partial_z, b(W) = W, \quad \forall \, W \in M_n(\Cset).  
\]
The algebras $K_1$ and $K_2$ are the algebras of matrix-valued polynomials in $x$ and $z$. The 
algebras $A_1$ and $A_2$ are the algebras of matrix differential operators in $x$ and $z$ with 
constant coefficients. The classification result is as follows:

\bth{th-i2} The bispectral Darboux transformations from the function $e^{xz}I_n$ as in \thref{th-i} 
for operators $L, Q, P$ with invertible leading terms are precisely the 
functions of the form
\[
\Phi(x,z) = W(x) P(x, \partial_x) (e^{xz} I_n),
\] 
where $W(x)$ is an invertible matrix-valued rational function and
$P(x, \partial_x)$ is a monic matrix differential operator of order $k$ that has a nondegenerate vector
kernel with a basis of the form
\[
e^{\al_1 x} p_1(x), \ldots, e^{\al_{kn} x} p_{kn}(x)
\]
for some (not necessarily distinct) $\al_j \in \Cset$ and vector valued polynomials $p_j(x) \in \Cset[x]^n$.
\eth

The operators $P(x, \partial_x)$ are uniquely reconstructed from their vector 
kernels using quasideterminants, see Sect. \ref{3.2}. The restriction to operators $L, Q, P$ with 
invertible leading terms is necessary because of the standard pathological problems with 
leading terms that are zero divisors. See \deref{non-deg-subsp} for the notion of nondegenerate space of vector-valued 
meromorphic functions.

In \cite{Wi1,Wi2}, all rank 1 scalar bispectral functions were shown to be parametrized by the points of the adelic Grassmannian and, equivalently, the points of 
all Calogero--Moser spaces. The classification in \thref{th-i2} and Theorem 4.2 in \cite{Wi3} can be used to show that the above class 
coincides with the recently constructed classes of matrix bispectral functions by Wilson \cite{Wi3} and Bergvelt, Gekhtman and Kasman \cite{BGK} in 
relation to certain solutions of the multicomponent KP hierarchy and the Gibbons--Hermsen system parametrized by 
the vector adelic Grassmannian. (More precisely, to match the two families, our bispectral functions should be multiplied by $W(x)^{-1}$ and by an invertible
matrix-valued rational function in $z$ to make the leading terms of their asymptotic expansions at $\infty$ equal to $e^{xz} I_n$.) 
One should note that, even though 
the starting function $\Psi(x,z)= e^{xz} I_n$ is scalar matrix-valued, the set of bispectral matrix Darboux transformations from 
it is much more complicated than the set of bispectral Darboux transformation from $e^{xz}$, the latter being 
the union of all Calogero--Moser spaces \cite{Wi2}.  
\subsection{}
\label{1.4}
Our second classification result is for all bispectral matrix Darboux trasnformations from 
the Airy functions. The (generalized) Airy operator is the differential operator 
\[
M_{\Ai}(x, \partial_x) := \partial^N_x + \sum_{i=1}^{N-1}\alpha_i \partial^{N-i}_x + \alpha_0 x,
\]
where $\alpha_0 \neq 0, \alpha_2, \ldots, \alpha_{N-1}$ are complex numbers. 
Its kernel is $N$-dimensional over $\Cset$. For any function $\psi(x)$ in the kernel of $M_{\Ai}(x, \partial_x)$, 
define the matrix Airy bispectral function 
\[
\Psi_{\Ai}(x,z) := \psi(x+z) I_n. 
\]
The bispectral Darboux transformations from it in the scalar case $(n=1)$ 
were classified in \cite{BHY2,KR} and played a key role in the treatment of the scalar bispectral problem.
Here we obtain a classification of these transformations for all $n$. 

In the Airy case, we again have $\Om_1 = \Om_2 = \Cset$. The algebras $B_1$ and $B_2$ are also  
the algebras of matrix-valued differential operators with polynomial coefficients in $x$ and $z$,
but the isomorphism $b \colon B_1 \to B_2$ is given by 
\[
b(x) := M_{\Ai}(z, \partial_z), \; \; b(\partial_x):= - \partial_z, \; \; b(W) := W, \; \forall \, W \in M_n(\Cset). 
\]
Here and below $M_{\Ai}(x,\partial_x)$ is viewed as a scalar matrix-valued differential operator. 
The first two equations imply that $b(M_{\Ai}(x, \partial_x)) = z$. The definition of the isomorphism $b$ 
represents the bispectral equations for the matrix Airy function
\[ 
M_{\Ai}(x, \partial_x) \Psi_{\Ai}(x,z) = \Psi_{\Ai}(x,z) z, \quad
x \Psi_{\Ai}(x,z) = \Psi_{\Ai}(x,z) M_{\Ai}(z, \partial_z).
\]
The algebras $K_1$ and $K_2$ associated to $\Psi_{\Ai}(x,z)$ 
are also the algebras of matrix-valued polynomial functions in $x$ and $z$, 
respectively. However, in contrast with \S \ref{1.3}, 
the algebras $A_1$ and $A_2$ consist of the operators of the form 
\[
q(M_{\Ai}(x, \partial_x)) \quad \mbox{and} \quad
q(M_{\Ai}(z, \partial_z))
\]
for matrix-valued polynomials $q(t)$. Denote by 
\[
\psi_1(x), \ldots, \psi_N(x)
\]
a basis of the kernel of $M_{\Ai}(x, \partial_x)$. We have:

\bth{th-i3} The bispectral Darboux transformations from the matrix Airy function $\Psi_{\Ai}(x,z)$ 
as in \thref{th-i} for operators $L, Q, P$ with invertible leading terms are exactly 
the functions of the form
\[
\Phi(x,z) = W(x) P(x, \partial_x) \Psi_{\Ai}(x,z),
\] 
where $W(x)$ is an invertible matrix-valued rational function and 
$P(x, \partial_x)$ is a monic matrix differential 
operator of order $dN$ whose vector kernel is a nondegenerate $ndN$-dimensional space 
with a basis which is a union of $N$-tuples of the form  
\[
\sum_{j=0}^{N-1}\psi_i^{(j)} (x + \la) p_j(x), \quad i=1, \ldots N
\]
for some matrix polynomials $p_j(x)$. (The matrix polynomials and complex numbers $\la \in \Cset$ 
can be different for different $N$-tuples of basis elements.) 
\eth

We expect that the families of bispectral Darboux transformations from the matrix Airy functions
will have applications to the study of the geometry of quiver varieties and to tau-functions 
of integrable systems satisfying the string equation.
\subsection{}
\label{1.5}
The paper is organized as follows. The proof of \thref{th-i} and its general form for bispectral Darboux 
transformations in bimodules of noncommutative algebras is given in Sect. \ref{GenTh}. This section also 
contains a review of the needed facts from noncommutative algebra. The classification 
results for matrix rank 1 Darboux transformations are in Sect. \ref{rank1}. The classification results for 
bispectral Darboux transformations from the matrix Airy function appear in Sect. \ref{Airy}. Sect. 
\ref{exe} contains examples illustrating the classification results. 

The classification results rely on some facts on quasideterminants and the spectral theory 
of matrix polynomials. The former are reviewed in Sect. \ref{Facto} and the latter in the appendix, 
Sect. \ref{Appen}.
\medskip
\\
\noindent
{\bf Acknowledgements.} We would like to thank Alberto Gr\"unbaum and Alex Kasman for their very helful
comments on the first version of the paper. M.Y. thanks the hospitality of Sofia University, Newcastle University and the 
Max Planck Institute for Mathematics in Bonn where parts of this project were completed. 
\sectionnew{General Theorem}
\label{GenTh}
In this section we establish a theorem that proves bispectrality of noncommutative Darboux transformations in an abstract algebraic setting. 
After the statement of the theorem we explain how it specializes to the concrete settings of continuous
and discrete matrix-valued bispectral operators.

\subsection{}
\label{2.1}
Before we proceed with the statement of the theorem we recall the notion of Ore set and noncommutative 
localization. We refer the reader to \cite[Ch. 6]{GW} for details.
Let $R$ be a unital (in general noncommutative) ring. A subset $E$ of $R$ is called multiplicative
if it is closed under multiplication and contains the identity element $1$. A multiplicative subset $E$ of $R$ 
consisting of regular elements (non-zero divisors) is called a (left) Ore set if for all $e \in E$ and $r \in R$
there exist $g \in E$ and $t \in R$ such that $g r = t e$.
If this condition is satisfied, then one can form the localized (or quotient) ring $R[E^{-1}]$ consisting 
of left fractions $e^{-1} r$, where $e \in E$, $r \in R$ with a certain natural equivalence relation. The ring $R$ 
canonically embeds in the localization by $r \mt 1^{-1} r$. 

For a ring $R$, we denote by $r(R)$ the multiplicative subset of its regular elements (i.e., elements that 
are not zero divisors).

\bth{GenTh}
Let $B_1$ and $B_2$ be two associative algebras and $M$ a $B_1-B_2$ bimodule. 
Assume that there exists an algebra isomorphism $b \colon B_1 \to B_2$ and 
an element $\psi \in M$ such that 
\[
P \psi = \psi b(P), \quad \forall P \in B_1.
\] 
Assume that $A_i$ and $K_i$ are subalgebras of $B_i$ such that 
\[
b(A_1) = K_2 \quad \mbox{and} \quad b(K_1) = A_2
\]
and for which the following conditions are satisfied: 
\begin{enumerate}
\item
$r(K_i)$ are Ore subsets of $B_i$ consisting 
of regular elements and

\item $\Ann_f(M) = 0$ for all $f \in r(K_i)$.  
\end{enumerate}
If $L = Q g^{-1} P$ for some $L \in A_1$, $P,Q \in B_1$, and $g \in r(K_1)$, then 
$\vph := P \psi = \psi b(P)$ satisfies
\begin{align*}
\left( P Q g^{-1} \right) \vph &= \vph b(L), \\
g \vph &= \vph \left( b(L)^{-1} b(Q) b(P) \right).
\end{align*} 
\eth
Note that $b(L) \in K_2$. The two conditions ensure that 
the algebras $B_i$ embed in $B_i[r(K_i)^{-1}]$ and 
$M$ in $M[r(K_1)^{-1}][r(K_2)^{-1}]$,
see e.g. \cite[Ch. 10]{GW} for details. Both equalities take place in 
the localized module. 

We prove the theorem in Sect. \ref{2.4}. In the next two subsections we 
describe how it is applied in the continuous--continuous and discrete-discrete 
noncommutative bispectral situations.

\subsection{}
\label{2.2}
In the continuous--continuous noncommutative bispectral setting the theorem works as follows.
Let $R$ be a complex finite dimensional algebra. For $i=1,2$, let $\Om_i$ be two open subsets of $\Cset$. 
Let 
\[
\Psi \colon \Om_1 \times \Om_2 \to R
\]
be an $R$-valued bispectral function. Recall from the introduction that 
this means that there exist $R$-valued analytic differential operators 
$L(x, \partial_x)$ and $\La(z, \partial_z)$ on 
$\Om_1$ and $\Om_2$, and two $R$-valued analytic functions 
\[
\theta \colon \Om_1 \to R, \quad 
f \colon \Om_2 \to R, 
\]
such that
\begin{align}
\label{bisp1}
L(x, \partial_x) \Psi(x,z) &= \Psi(x,z) f(z) \\
\label{bisp2}
\theta(x) \Psi(x, z) &=
\Psi(x,z) \La(z, \partial_z)
\end{align}
on $\Om_1 \times \Om_2$. We will furthermore assume that $\Psi(x,z)$ is a nonsplit function of $x$ and $z$
in the sense that it satisfies the condition

$(*)$ there are no nonzero analytic differential operators $L(x, \partial_x)$ and $\La(z, \partial_z)$
that satisfy one of the above conditions with $f(z) \equiv 0$ or $\theta(x) \equiv 0$.

In this setting, denote by $B_1$ the algebra of all $R$-valued analytic differential operators $P(x, \partial_x)$ 
on $\Om_1$ for which there exits an $R$-valued analytic differential operator $S(z, \partial_z)$ on $\Om_2$ such that 
\[
P(x, \partial_x) \Psi(x,z) = \Psi(x,z) S(z, \partial_z).
\]
Denote by $B_2$ the algebra of all such operators $S(z, \partial_z)$. It follows from the assumption $(*)$
that the map
$b \colon B_1 \to B_2$, given by $b(P(x, \partial_x)):= S(z, \partial_z)$, is 
a well defined algebra isomorphism. Let $M$ be the $B_1 - B_2$ bimodule consisting of 
$R$-valued analytic functions on $\Om_1 \times \Om_2$.

Let $K_i$ be the subalgebras of $B_i$ consisting of $R$-valued analytic functions. Set
\[
A_1 := b^{-1}(K_2), \quad A_2 := b(K_1).
\]
These are the algebras of bispectral operators corresponding to $\Psi(x,z)$. They consist 
of the differential operators that have the properties 
\eqref{bisp1}--\eqref{bisp2} for some functions $f(z)$, $\theta(x)$.  
In particular,  
\[
L(x, \partial_x) \in A_1 \quad \mbox{and} \quad 
\La(z, \partial_z) \in A_2.
\]
The sets of regular elements $r(K_i)$ consist of those $R$-valued analytic functions $\theta(x)$ and $f(z)$ 
that are nondegenerate (in the sense that their determinants 
are not identically zero on any connected component of $\Om_i$).
This immediately implies that the condition (2) in \thref{GenTh}
is satisfied in this setting. By standard commutation arguments with differential operators one shows 
that the subsets $r(K_i)$ are Ore subsets of $B_i$.   

\thref{GenTh} implies the following:

\bpr{p1} In the above setting, for every factorization, 
\[
L(x, \partial_x) = Q(x, \partial_x) g(x)^{-1} P(x, \partial_x)
\]
for $L(x, \partial_x) \in A_1$, $P(x, \partial_x), Q(x, \partial_x) \in B_1$, $g(x) \in r(K_1)$, 
the $R$-valued function 
\[
\Phi(x,z) := P(x,\partial_x) \Psi(x,z) \colon \Om_1 \times \Om_2 \to R
\]
is bispectral, and, more precisely, satisfies
\begin{align*}
\left( P(x, \partial_x) Q(x, \partial_x)  g(x)^{-1} \right) \Phi(x,z) 
&= \Phi(x,z) b(L)(z) \\
g(x) \Phi(x,z) &= \Phi(x,z) \left( b(L)(z)^{-1} b(Q)(z, \partial_z) b(P)(z, \partial_z) \right).
\end{align*} 
\epr 

\subsection{}
\label{2.3}
In the discrete--discrete noncommutative bispectral setting \thref{GenTh} works as follows. As before,
let $R$ be a complex finite dimensional algebra. 
For $i=1,2$, let $\Om_i$ be two subsets of $\Cset$ which are invariant 
under the translation operator
\[
T_x \colon x \mt x+1
\]
and its inverse.
Denote by $M$ the space of $R$-valued functions 
on $\Om_1 \times \Om_2$. An $R$-valued difference operator on $\Om_1$ 
is a finite sum of the form
\[
\sum_{k \in \Zset} c_k(x) T_x^k, 
\]
where $c_k$ are $R$-valued functions. The space 
$M$ is naturally equipped with the structure of a bimodule
over the algebras of $R$-valued difference operators 
on $\Om_1$ and $\Om_2$, respectively, by setting
\begin{align*}
\left( \sum_{k \in \Zset} c_k(x) T_x^k \right) \cdot \Psi(x, z) &:= 
\sum_{k \in \Zset} c_k(x)  \Psi(x+k, z), 
\\
\Psi(x,z) \cdot \left( \sum_{k \in \Zset} c_k(z) T_z^k \right) &:= 
\sum_{k \in \Zset} c_k(z-k)  \Psi(x, z-k).
\end{align*}
An $R$-valued discrete-discrete bispectral function is by definition an element of M (i.e., 
an $R$-valued function) 
\[
\Psi \colon \Om_1 \times \Om_2 \to R
\]
for which there exist $R$-valued difference operators 
$L(x, T_x)$ and $\La(z, T_z)$ on 
$\Om_1$ and $\Om_2$, and $R$-valued functions 
\[
\theta \colon \Om_1 \to R, \quad 
f \colon \Om_2 \to R, 
\]
such that
\begin{align}
\label{bisp3}
L(x, T_x) \Psi(x,z) &= \Psi(x,z) f(z) \\
\label{bisp4}
\theta(x) \Psi(x, z) &=
\Psi(x,z) \La(z, T_z)
\end{align}
on $\Om_1 \times \Om_2$. Similarly to the continuous case, we will assume that $\Psi(x,z)$ is a nonsplit 
function of $x$ and $z$ in the sense that it satisfies the condition

$(**)$ there are no nonzero difference operators $L(x, T_x)$ and $\La(z, T_z)$
that satisfy one of the above conditions with $f(z) \equiv 0$ or $\theta(x) \equiv 0$.

Let $B_1$ be the algebra of all $R$-valued difference operators $P(x, T_x)$ 
on $\Om_1$ for which there exits an $R$-valued difference operator $S(z, T_z)$ on $\Om_2$ such that 
\[
P(x, T_x) \Psi(x,z) = \Psi(x,z) S(z, T_z).
\]
Let $B_2$ be the algebra of all such operators $S(z, \partial_z)$. The assumption $(**)$
implies that the map $b \colon B_1 \to B_2$, given by $b(P(x, \partial_x)):= S(z, \partial_z)$ is 
a well defined algebra isomorphism. Denote by $K_i$ be the subalgebras of $B_i$ consisting of $R$-valued 
functions. The algebras
\[
A_1 := b^{-1}(K_2), \quad A_2 := b(K_1)
\]
consist of the bispectral difference operators corresponding to $\Psi(x,z)$ (i.e., 
difference operators in $x$ and $z$ that have the properties \eqref{bisp3}--\eqref{bisp4}). 
The sets $r(K_i)$ of regular elements of the algebras $K_i$ consist of the nonvanishing 
functions $\theta(x)$ and $f(z)$. It is easy to verify that the second condition in 
\thref{GenTh} is satisfied in this situation and that 
$r(K_i)$ are Ore subsets of $B_i$. We leave the details to the reader.

\thref{GenTh} implies the following:

\bpr{p2} In the above setting, for every factorization
\[
L(x, T_x) = Q(x, T_x) g(x)^{-1} P(x, T_x)
\]
with $L(x, T_x) \in A_1$, $P(x, T_x), Q(x, T_x) \in B_1$, $g(x) \in r(K_1)$, 
the function 
\[
\Phi(x,z) := P(x, T_x) \Psi(x,z)
\]
is an $R$-valued discrete-discrete bispectral function on $\Om_1 \times \Om_2$. 
More precisely, it satisfies
\begin{align*}
\left( P(x, T_x) Q(x, T_x)  g(x)^{-1} \right) \Phi(x,z) 
&= \Phi(x,z) b(L)(z), \\
g(x) \Phi(x,z) &= \Phi(x,z) \left( b(L)(z)^{-1} b(Q)(z, T_z) b(P)(z, T_z) \right).
\end{align*} 
\epr

Analogously, \thref{GenTh} applies to noncommutative bispectral functions for $q$-difference 
operators and to every mixed situation, where the operators on the two sides 
are differential, difference, or $q$-difference of different type. We leave to the reader 
the formulations of those results, which are analogous to the ones above.

The bispectral functions in \cite{BL,G1,Ka,Z1} are special cases of the continuous-continuous version 
of the theorem and the ones in   \cite{Du, DuG, GI, GPT, GPT2, M}  are special cases of the discrete-continuous version. 

\subsection{}
\label{2.4}
We finish this section with the proof of \thref{GenTh}.
\medskip
\\
\noindent
{\em{Proof of \thref{GenTh}.}}
The first statement is straightforward, 
\[
\left( P Q g^{-1} \right) \vph = 
\left( P Q g^{-1} \right) P \psi = 
P L \psi = P \psi b(L) = \vph b(L).
\]

The key point of the theorem is the second statement. To prove it, first 
note that since $b \colon B_1 \to B_2$ is an isomorphism, the subalgebras 
$A_i$ of $B_i$ also satisfy the conditions (1)--(2). The Ore subset assumptions 
imply that
\[
g = P L^{-1} Q \quad \mbox{and thus}\quad
b(g) = b(P) b(L)^{-1} b(Q).
\]
Therefore,
\[
g \vph = g \psi b(P) = \psi b(P) b(L)^{-1} b(Q) b(P) = \vph \left( b(L)^{-1} b(Q) b(P) \right).
\]
\qed
\sectionnew{Quasideterminants, kernels and factorizations of matrix differential operators}
\label{Facto}
In this section we gather some facts relating the factorizations of matrix-valued differential operators 
and their kernels, which will be needed in the following sections. These results are minor 
variations of those obtained by Etingof, Gelfand and Retakh in \cite{EGR}. One difference in our treatment 
is that we focus on vector kernels rather than matrix kernels (considered in \cite{EGR}) because 
our classification results for bispectral Darboux transformations from the next 2 sections are 
naturally formulated in terms of the former.

\subsection{}
\label{3.1}
We first recall the notion of quasideterminants introduced and studied by Gelfand and Retakh 
\cite{GR1,GR2}.
Let $R$ be an associative algebra. Let $X=(x_{ij})$ be a matrix with elements in $R$. For each pair of indices $i, j$ denote by $r_i(X)$ the $i$-th row of $X$ and by $c_j(X)$ the $j$-th column of $X$. Denote by $X^{ij}$ the matrix obtained by removing the $i$-th row and the $j$-th column of $X$. For each row $r$ and index $j$, define the vector $r^{(j)}$ obtained by removing the $j$-th entry. Similarly, for a column $c$, denote by $c^{(i)}$ the column vector obtained by its $i$-th entry. Finally, we assume that the  matrix $X^{ij}$ is invertible. The quasideterminant  $|X|_{ij}$ 
is defined by 
\[
|X|_{ij}   :=    x_{ij} - r_i(X)^{(j)} (X^{ij})^{-1} c_j(X)^{(i)} \in R.
\]

Let $\Om \subseteq \Cset$ be a domain and $R$ be the algebra of meromorphic functions on
$\Om$ with values in $M_n(\Cset)$. 

For any sets of matrix-valued meromorphic
functions $F_1, \ldots, F_k$ on $\Om$ and 
vector-valued meromorphic
functions $f_1, \ldots, f_{nk}$ on $\Om$,
define the Wronski matrices 
\begin{equation}
\label{Wr}
W(F_1, \ldots, F_k) =   \left(
\begin{matrix}
F_1 & \ldots & F_k \\ 
F'_1 & \ldots & F'_k  \\ 
\ldots & \ldots & \ldots \\ 
F^{(k-1)}_1 & \ldots &F^{(k-1)}_k
\end{matrix} \right), 
W(f_1, \ldots, f_{nk}) =   \left(
\begin{matrix}
f_1 & \ldots & f_{nk} \\ 
f'_1 & \ldots & f'_{nk}  \\ 
\ldots & \ldots & \ldots \\ 
f^{(nk-1)}_1 & \ldots & f^{(nk-1)}_{nk}
\end{matrix} \right). 
\end{equation}
Such a set of functions $F_1, \ldots, F_k$ is called {\em{nondegerate}} if $W(F_1, \ldots, F_k)$ is 
an invertible matrix-valued meromorphic function on $\Om$.
The invertibility condition is independent of whether \eqref{Wr} is considered as an $n \times n$ matrix with values 
in $R$ or an $nk \times nk$ matrix whose entries are meromorphic functions on $\Om$. This is equivalent to saying that 
the determinant of the Wronski $nk \times nk$ matrix is a nonzero meromorphic function on $\Om$.

\bde{non-deg-subsp} A subspace of dimension $nk$ of the space of vector-valued meromorphic functions on $\Om$ will be 
called {\em{nondegenerate}} if it has a basis $f_1, \ldots, f_{nk}$ such that 
the Wronski matrix $W(f_1, \ldots, f_{nk})$ is an invertible matrix-valued meromorphic function on $\Om$. 
\ede
It is clear that the condition in the definition does not depend on the choice of basis. 
\bre{arrange}
For every basis 
of a nondegenerate subspace of the space of vector-valued meromorphic functions on $\Om$, 
there exists an arrangement of the basis elements into the columns of 
a set of meromorphic functions $F_1, \ldots, F_k \colon \Om \to M_n(\Cset)$ with the property that 
for all $m \leq k$ the set $F_1, \ldots, F_m$ is nondegenerate on $\Om$. This follows from the fact 
that, for every nondegenerate square matrix $M$, there exists a permutation matrix $\Sig$ such that all leading principal minors 
of $M \Sig$ are nondegenerate (the fact is applied for $M$ equal to the Wronski matrix of the basis elements).
\ere
\subsection{}
\label{3.2}
Define the vector kernel of a differential operator $P(x,\partial_x)$ with matrix-valued meromorphic
coefficients on $\Om$ by
\begin{equation}
\label{vker}
\vker P(x,\partial_x) := \{ \mbox{meromorphic functions} \; \; f(x) \colon \Om \to \Cset^n \mid P(x, \partial_x) f(x) =0 \}.
\end{equation}
By the standard theorem for existence and uniqueness of solutions of ordinary differential equations, 
$\dim \vker P(x, \partial_x) \leq kn$, and each $a \in \Om$ at which $P(x, \partial_x)$ is regular 
has a neighborhood $\OO_a \subseteq \Om$
such that the restriction of $P(x, \partial_x)$ to $\OO_a$ satisfies 
$\dim \vker P(x, \partial_x) = kn$.

\bpr{diff-op} Let $\Om \subseteq \Cset$ be a domain and $n, k$ be positive integers.

(a) Assume that $P(x, \partial_x)$ is an $M_n(\Cset)$-valued 
meromorphic differential operator on $\Om$ of order $k$ that has $nk$-dimensional vector kernel.
Then its kernel is nondegenerate in the sense of \deref{non-deg-subsp}. 

(b) Assume that $V$ is a nondegenerate subspace of the space of vector-valued meromorphic
functions on $\Om$ of dimension $nk$.  Then there is a unique monic differential operator $P(x, \partial_x)$ of order $k$ with matrix-valued meromorphic 
coefficients on $\Om$ whose vector kernel equals $V$. For every basis of $V$ choose an arrangement of the basis elements onto the columns of 
a set of meromorphic functions $F_1, \ldots, F_k \colon \Om \to M_n(\Cset)$ with the property that 
for all $m \leq k$ the set $F_1, \ldots, F_m$ is nondegenerate on $\Om$ (cf. \reref{arrange}). The operator $P(x, \partial_x)$ is given by 
\begin{equation}
\label{Pf}
P(x, \partial_x) F(x) = |W(F_1, \ldots, F_k, F)|_{k+1, k+1}
\end{equation}
for each matrix-valued meromorphic function $F(x)$ and by
\begin{equation}
P(x, \partial_x) = (\partial_x - b_k) \ldots (\partial_x - b_1), 
\label{L-form2}
\end{equation}
where
\begin{equation}
\label{bj}
b_j := W'_j W_j^{-1}, \; \; W_j := |W(F_1, \ldots, F_j)|_{jj}. 
\end{equation}
\epr
\begin{proof} (a) Denote by $\Om^\circ$ the open dense subset of $\Om$ on which $P(x, \partial_x)$ is analytic. 
Choose a basis $f_1(x), \ldots, f_{nk}(x)$ of $\vker P(x, \partial_x)$. Let  $c_0, \ldots, c_{nk} \in \Cset$. 
The statement of this part of the proposition follows from the fact that if 
\[
W(f_1, \ldots, f_{nk})(x_0) \cdot (c_1, \ldots c_{nk})^t = 0  
\]
for some $x_0 \in \Om^\circ$, then this holds for all $x_0 \in \Om^\circ$, where $(.)^t$ denotes the transpose of a matrix.
(This is nothing 
but the uniqueness statement for solutions of the initial value problem for ordinary differential 
equations.) 

(b) By \cite[Theorem 1.1]{EGR} there exists a unique monic differential operator $P(x, \partial_x)$ of order $k$ such that 
\[
P(x, \partial_x) F_i(x) = 0, \quad \forall \; 1 \leq i \leq k
\]
and this operator is given by \eqref{Pf} and \eqref{L-form2}. Since $\vker P  \supseteq V$,
$nk \geq \dim \vker P$ and $\dim V = n k$, we have  $\vker P =V$. 
\end{proof}
\subsection{} The factorizations of matrix-valued differential operators are described by the next proposition.
\label{3.3}
\bpr{fact} Let $L(x, \partial_x)$ be an $M_n(\Cset)$-valued monic meromorphic differential operator on a 
domain $\Om \subseteq \Cset$ of order $l$ whose vector kernel has dimension $nl$. All factorizations 
\begin{equation}
\label{facto}
L(x, \partial_x) = Q(x, \partial_x) P(x, \partial_x) 
\end{equation}
into $M_n(\Cset)$-valued monic meromorphic differential operators of orders $l-k$ and $k$ 
are in bijection with the nondegenerate $nk$-dimensional subspaces of $vker L$. 
This bijection is given by $P \mt \vker P$. The inverse bijection from nondegenerate subspaces of $\vker L$ 
to monic differential operators $P$ satisfying \eqref{facto} is given by \prref{diff-op} (b).
\epr  
This proposition is similar to \cite[Theorem 1]{GV}, but the latter omits the condition for nondegeneracy of $\vker P$. 

\begin{proof} Assume that we have a factorization as in \eqref{facto}. Obviously, $\vker P \subseteq \vker L$. 
We will show that $\dim \vker P = nk$. Then by \prref{diff-op} (a), $\vker P$ is nondegenerate.  

Fix a point $a \in \Om$ at which $L(x,\partial_x)$ 
is regular and a neighborhood 
$\OO_a$ of it on which $Q$ and $P$ are analytic and the restriction of $P$ has $nk$-dimensional kernel.
For this restriction $\vker P \subseteq \vker L$, i.e., there exists an $nk$-dimensional 
subspace $V$ of the vector kernel of $L(x, \partial_x)$ such that the vector kernel of the restriction of $P$ 
to $\OO_a$ is $V|_{\OO_a}$. By applying \prref{diff-op} (b), we construct a differential operator $P'(x, \partial_x)$ on $\Om$ with meromorphic coefficients
such that $\vker P' = V$.  The restrictions of $P$ and $P'$ to $\OO_a$ coincide by \prref{diff-op} (b) and 
the two differential operators are meromorphic on $\Om$. Since $\Om$ is connected,
$P = P'$. So, $\vker P = V$, and, in particular, $\dim \vker P = nk$.

In the other direction, for an $nk$-dimensional subspace $V$ of the vector kernel of $L$ on $\Om$ and a basis of it, 
we can form meromorphic functions $F_1, \ldots, F_k \colon \Om \to M_n(\Cset)$ from a basis of $V$ as in \reref{arrange}.
Furthermore, one can complement 
the basis of $V$ to a basis of $vker L$ and form a set of meromorphic functions 
$F_{k+1}, \ldots, F_l \colon \Om \to M_n(\Cset)$ from them, such that 
$F_1, \ldots F_m$ is a nondegenerate set for all $m \leq l$.  
Then, by \eqref{L-form2},
\[
L(x, \partial_x) = Q(x, \partial_x) P(x, \partial_x) 
\]  
for 
\[
Q(x, \partial_x) := (\partial_x - b_l) \ldots (\partial_x - b_{k+1}),\quad
P(x, \partial_x) := (\partial_x - b_k) \ldots (\partial_x - b_1), 
\]
where $b_j$ are given by \eqref{bj}.
\end{proof}
In the setting of the above proposition the meromorphic operator $PQ$ is called a 
matrix Darboux transformation from $L = QP$. For an eigenfuction $\Psi(x,z)$ of $L(x, \partial_x)$
(satisfying $L(x, \partial_x) \Psi(x,z) = \Psi(x,z) f(z)$) 
the function $P(x, \partial_x) \Psi(x,z)$ is also called a
Darboux transformation from $\Psi(x,z)$. 
 
\bre{LP} Some authors consider the following form of Darboux transformations. 
For a given matrix-valued monic meromorphic 
differential operator $L(x, \partial_x)$, consider all matrix-valued
meromorphic operators $\ol{L}(x, \partial_x)$ that have the property 
\begin{equation}
\label{PL}
P(x, \partial_x) L(x, \partial_x) = \ol{L}(x, \partial) P(x, \partial_x)
\end{equation}
for some matrix-valued monic meromorphic operators $P(x, \partial_x)$. As before, we call
the function $P(x, \partial_x) \Psi(x,z)$ a 
Darboux transformation from $\Psi(x,z)$.

A priori it appears that this form of Darboux transformations is more general than the ones we consider.
However, this is not the case, i.e., \thref{GenTh} also covers these transformations.
To see this, one first proves that \eqref{PL} implies $L (\vker P) \subseteq \vker P$. Since 
$\dim \vker P < \infty$, $\vker P \subseteq \vker h(L)$ for some 
polynomial $h(z)$. Therefore $h(L) = Q P$ for some matrix-valued meromorphic differential 
operator $Q$ on $\Om$ and we have the Darboux transformation 
\[
h(L) = QP \mt PQ = h(\ol{L})
\]
which falls within the class treated in \thref{GenTh}.
On the level of wave functions it represents the Darboux transformation \eqref{PL}
\[
\Psi(x,z) \mt P(x, \partial_x) \Psi(x,z).
\]
\ere
\sectionnew{Classification of matrix bispectral Darboux transformations in the rank 1 case}
\label{rank1}
All bispectral algebras of rank 1 scalar ordinary differential operators were classified 
by Wilson \cite{Wi1}. In \cite{BHY2} it was proved that the corresponding wave functions
are bispectral Darboux transformations from the function $\Psi(x,z) = e^{xz}$ in the sense 
of \thref{GenTh}. At the time the result was phrased in different terms, but it is easy 
to convert the statements of \cite{BHY2} to the ones in \thref{GenTh}, and more precisely,
to the case of \prref{p1} when the size of the matrices is $1 \times 1$. 

In this section we classify explicitly all bispectral Darboux transformations 
of the corresponding matrix analog in full generality, i.e., we obtain a full classification 
of all matrix rank 1 bispectral Darboux transformations. This is achieved in \thref{class-rank1}.
The only restriction we impose is that the leading terms of the differential operators $L, Q$ and $P$ in \prref{p1} 
be invertible matrix-valued meromorphic functions. Without loss of generality this situation is the same as the one
when $L, Q$ and $P$ are monic. (The classifications of the factorizations of differential operators with 
noninvertible leading terms are related to well known pathological problems from the zero divisors of the algebra 
$M_n(\Cset)$ and will be considered elsewhere.) 
\subsection{}
\label{4.1}
In the setting of the previous section, let $\Om_1 = \Om_2 := \Cset$. Define the function
\[
\Psi(x,z) = e^{xz} I_n
\] 
on $\Cset \times \Cset$, where $I_n$ is the identity matrix. Let $B_1$ and $B_2$ be both 
equal to the algebra of matrix-valued differential operators with polynomial coefficients. 
We will use the variables $x$ and $z$ for the first and second algebra, respectively.
Let $b \colon B_1 \to B_2$ be the isomorphism given by 
\[
b(\partial_x) := z, \; \; b(x) = - \partial_z, \; \; b(W) := W, \; \forall \, W \in M_n(\Cset). 
\]
We are in the bispectral situation of \prref{p1}:
\[
P(x, \partial_x) \Psi(x,z) = \Psi(x,z) (bP)(z, \partial_z), \; \; \forall P(x, \partial_x) \in B_1. 
\]
Let $K_1$ and $K_2$ be the algebras of matrix-valued polynomial functions and $x$ and $z$ respectively. 
Denote by $A_1$ and $A_2$ the algebras of matrix-valued differential operators (in $x$ and $z$) with 
constant coefficients. Then 
\[
b(A_1) = K_2, \; \; b(K_2) = A_1.
\]
The regular elements of the algebras $K_1$ and $K_2$ are the matrix-valued 
polynomial functions $\theta(x)$ and $f(z)$ that are nondegenerate in the sense that 
their determinants are nonzero polynomials.
By standard commutation arguments for matrix differential operators the sets 
$r(K_i)$ are Ore subsets of $B_i$. 
The localizations $B_1[r(K_1)^{-1}]$ and $B_2[r(K_2)^{-1}]$
are nothing but the algebras of differential operators with rational coefficients 
in the variables $x$ and $z$. \thref{GenTh} implies the following:

\bth{rank1} If $L(x, \partial_x)$ is a matrix-valued differential operator with constant coefficients, 
and $Q(x, \partial_x)$ and $P(x, \partial_x)$ are matrix-valued differential operators with rational coefficients, 
such that
\begin{equation}
\label{rank1-fact}
L(x, \partial_x) = Q(x, \partial_x) P(x, \partial_x),
\end{equation}
then 
\[
\Phi(x,z) := P(x, \partial_x) (e^{xz} I_n) 
\]
is a matrix bispectral function.

More precisely, if 
\begin{equation}
\label{special-fact}
L(x, \partial_x) = Q'(x, \partial_x) g(x)^{-1} P'(x, \partial_x)
\end{equation}
for matrix-valued differential operators with polynomial coefficients $Q'(x, \partial_x)$, $P'(x, \partial_x)$ 
and a nondegenerate matrix-valued polynomial function $g(x)$,
then $\Phi(x,z) := P'(x, \partial_x) (e^{xz} I_n)$ satisfies
\begin{align*}
\left( P'(x, \partial_x) Q'(x, \partial_x) g(x)^{-1} \right) \Phi(x,z) &= \Phi(x,z) (b L)(z),
\\
g(x) \Phi(x,z) &= \Phi(x,z) \left( (bL)(z)^{-1} (bQ')(z, \partial_z) (bP')(z, \partial_z) \right).
\end{align*}
\eth

\subsection{}
\label{4.2}
Our goal below is to classify all wave functions $\Phi(x,z)$ that can be obtained from \thref{rank1} for operators $L(x, \partial_x)$
with nondegenerate leading term. This is done in \thref{class-rank1}. 
Denote the vector space of vector-valued  quasipolynomial functions
\begin{equation}
\label{QPn}
\QP_n := \oplus_{\al \in \Cset} e^{\al x} \Cset[x]^n,
\end{equation}
where $\Cset[x]$ denotes the space of polynomials and $\Cset[x]^n$ the space
of vector-valued polynomial functions.
The elements of $\QP_n$ have the form
$\sum_{\al} e^{\al x} p_{\al, i}(x)$, where the sum is finite and $p_{\al, i}(x) \in \Cset[x]^n$.

By converting $n$-th order vector differential equations to first order ones and using 
\prref{fact}, one easily obtains the following lemma. We leave its proof to the reader.
 
\ble{lin} (a) Assume that $L(x, \partial_x)$ is a monic matrix-valued differential operator with constant coefficients.
If $Q(x, \partial_x)$ and $P(x, \partial_x)$ are monic matrix-valued differential operators 
with meromorphic coefficients in $\Cset$ which satisfy \eqref{rank1-fact},  
then $\vker P(x, \partial_x)$ is a nondegenare subspace of $\QP_n$.

(b) For every nondegenerate finite dimensional subspace $V$ of $\QP_n$ whose dimension is a multiple of $n$, 
there exist monic matrix-valued differential 
operators $L$, $Q$ and $P$ satisfying the properties in part (a) 
such that $\vker P(x, \partial_x) = V$. 
\ele
 
The next theorem is our classification result of the matrix rank 1 bispectral Darboux transformations.
\bth{class-rank1} Assume that $L(x, \partial_x)$ is a monic matrix-valued differential operator with constant coefficients, 
and $Q(x, \partial_x)$ and $P(x, \partial_x)$ are monic matrix-valued differential operators whose coefficients are rational
functions. If \eqref{rank1-fact} is satisfied, then $\vker P(x, \partial_x)$ is a nondegenerate subspace of $\QP_n$ that has a basis of the form
\begin{equation}
\label{ker-P}
\{e^{\al_j x} p_j(x) \mid 1 \leq j \leq k n \},  
\end{equation}
where $\al_1, \ldots, \al_{kn}$ are not necessarily distinct complex numbers and $p_j(x) \in \Cset[x]^n$. 
Furthermore, for every such subspace of $\QP_n$ there exists a triple of operators 
$L$, $Q$ and $P$ with the above properties such that $v\ker P(x, \partial_x)$ equals this subspace.
\eth
\bco{class-rank1-b} The bispectral Darboux transformations from the function $e^{xz}I_n$ as in \thref{rank1} 
for operators $L, Q, P$ with nondegenerate leading terms are precisely the functions of the form
\[
\Phi(x,z) = W(x) P(x, \partial_x) (e^{xz} I_n),
\] 
where $P(x, \partial_x)$ is a monic matrix differential operator of order $k$ with vector kernel 
having a basis of the form \eqref{ker-P} and $W(x)$ is a matrix-valued nondegenerate rational function.
The operator $P(x, \partial_x)$ is given by \eqref{Pf}, 
where $F_1(x), \ldots, F_k(x)$ are the matrix-valued functions with columns 
\[
e^{\al_1 x} p_1(x), \ldots, e^{\al_{kn} x} p_{kn}(x),
\]
$\al_j \in \Cset$, $p_j(x) \in \Cset[x]^n$. 
\eco

The bispectral functions in this family were related to wave functions of solutions of the $n$-KP 
hierarchy by Wilson \cite{Wi3} and Bergvelt, Gekhtman, and Kasman \cite{Ka}. The fact that the family in \coref{class-rank1-b} is the same 
as the one in \cite{BGK,Wi3} (up to a multiplication on the left by $W(x)^{-1}$ and by an invertible matrix-valued rational function in $z$)
follows from Theorem 4.2 in \cite{Wi3} and \coref{class-rank1-b}. 
Since the proof of 
Theorem 4.2 in \cite{Wi3} is only sketched there, we note that the reference to it is not used in proofs here, but 
only to relate the bispectral functions in this section to those in \cite{BGK,Wi3}.
The bispectrality of the special elements of this family obtained 
by transformations from $L(x, \partial_x) = q(\partial_x) D$ for a diagonal matrix $D$ and 
a (scalar) polynomial $q(t)$ was previously obtained in \cite{BL}. Applications of \thref{GenTh} to
generalized bispectrality in relation to the results of Kasman \cite{Ka} are described in \S \ref{4.4}.
\subsection{}
\label{4.3}
We will break the proof of the \thref{class-rank1} into two lemmas. 
\ble{rational} Let, for $1 \leq j \leq kn$, $\al_j \in \Cset$ and $p_j(x)$ be vector-valued 
rational functions for which 
\begin{equation}
\label{ker-ap}
e^{\al_1 x} p_1(x), \ldots, e^{\al_{kn} x} p_{kn}(x)
\end{equation}
form a nondegenerate (and thus, linearly independent) set. Then the coefficients of the unique monic matrix-valued differential
operator with kernel spanned by \eqref{ker-ap} are rational functions.  
\ele
\begin{proof} We argue by induction on $k$. In the case $k=1$ the differential operator equals 
\[
\partial_x - F'(x) F^{-1}(x),
\]
where $F(x)$ is the matrix function with columns $e^{\al_1 x} p_1(x), \ldots, e^{\al_{n} x} p_{n}(x)$,
and (e.g. by Cramer's rule) $F'(x) F^{-1}(x)$ has rational coefficients.

In the general case, by \eqref{L-form2}, the differential operator has the form  
\[
P(x,\partial_x) =\ol{P}(x, \partial_x)(\partial_x - F'(x) F^{-1}(x)) 
\]
with $F(x)$ being nondegenerate and with columns as in \eqref{ker-ap}. Therefore, $\partial_x - F'(x) F^{-1}(x)$ has rational
coefficients and the operator $\ol{P}(x, \partial_x)$ has vector kernel spanned by 
\[
(\partial_x - F'(x) F^{-1}(x)) (e^{\al_j x} p_j(x)), \quad n < j \leq kn. 
\]  
By the inductive assumption, applied to $\ol{P}(x, \partial_x)$, the operator 
$\ol{P}(x, \partial_x)$ has rational coefficients, which implies the same property for $P(x,\partial_x)$.  
\end{proof}
\ble{ker-exp} Assume that $P$ is matrix-valued differential operator with coefficients that are rational functions. If 
\[
\sum_{s=1}^r e^{\lambda_s x}p_s(x) \in \vker P  
\]
for some vector-valued polynomial functions $p_s(x)$ and distinct complex numbers 
$\la_1, \ldots, \la_r$, then
\[
e^{\lambda_s x}p_s(x) \in \vker P, \quad \forall 1 \leq s \leq r. 
\]
\ele
\begin{proof} We have 
\[
\sum_{s=1}^r e^{\lambda_s x} \left( e^{-\lambda_s x} P(x, \partial_x) (e^{ \lambda_s x} p_s(x)) \right) = 0.
\]
Since the functions $e^{-\lambda_s x} P(x, \partial_x) (e^{ \lambda_s x} p_s(x))$ are rational, the nonzero 
functions in the sum are linearly independent because $\la_1, \la_2, \ldots, \la_r$ are distinct.
Thus all terms in the sum vanish.
\end{proof}
\noindent
{\em{Proof of \thref{class-rank1}.}} The first part of the theorem follows from Lemmas \ref{lker-exp} and \ref{llin} (a). 
The second part of the theorem follows from Lemmas \ref{lrational} and \ref{llin} (b).
\qed 
\subsection{}
\label{4.4} There are generalized versions of the differential bispectral problem from \eqref{eq3}--\eqref{eq4}. In one of them, one 
keeps the differential eigenvalue problem in \eqref{eq3} but replaces the operator $\La(z,\partial_z)$ in the right hand side of 
\eqref{eq4} with a mixed differential-translation operator. In this subsection, we briefly describe an application of 
\thref{GenTh} that gives a second proof of the results in \cite{Ka} on the construction of such generalized bispectral functions.

Firstly, for an analytic function $\Psi \colon \Om_1 \times \Om_2 \to M_n(\Cset)$ with an expansion
\[
\Psi(x,z) = \sum_{j,k=0}^\infty (x-x_0)^j a_{j,k} (z-z_0)^k, \quad a_{j,k} \in M_n(\Cset), x_0 \in \Om_1, z_0 \in \Om_2
\]
and a matrix $W \in M_n(\Cset)$, define the translation action
\[
\Psi(x,z) \cdot (T_z^W) := \sum_{j,k=0}^\infty (x-x_0)^j a_{j,k} (z- W -z_0)^k
\]
whenever the series converges. A mixed differential-translation operator is an operator obtained 
by composing differential operators with translation operators of this kind.

For an invertible matrix $H \in M_n(\Cset)$, denote by $Z(H) =\{ W \in M_n(\Cset) \mid WH = H W \}$ its centralizer.
Consider the analytic function
\[
\Psi_H(x, z) := e^{xz H}.
\]
Let $B_1$ be the algebra of differential operators in $x$ with coefficients that are quasipolynomial 
functions with values in $Z(H)$ (i.e., functions with values in $Z(H)$ whose matrix entries are finite sums of terms of 
the form $x^j e^{\al x}$, $j \in \Nset$, $\al \in \Cset$).
Let $B_2$ be the algebra consisting of mixed differential-translation operators of the form
\[
\sum_{j \in \Nset, c \in \Cset} a_{j,c}(z) \partial_z^j T_z^{c H} 
\]
given by finite sums and terms $a_{j,c}(z)$ that are polynomials with values in $Z(H)$. The following define
an isomorphism $b \colon B_1 \to B_2$ 
\[
b(x) = - \partial_z, \; b(\partial_x) = z, \; b(e^{ax}) = T_z^{-aH^{-1}}, \; b(W) = W, \; \; \forall W \in Z(H)
\]
and the function $\Psi_H(x,z)$ satisfies
\begin{equation}
\label{bH}
P \cdot \Psi_H(x,z) = \Psi_H \cdot b(P), \quad \forall P \in B_1.
\end{equation}
Let $A_1$ be the subalgebra of $B_1$ consisting of constant differential operators in $x$ with values in $Z(H)$ and $K_2$  be the 
subalgebra of $B_2$ consisting of polynomial in $z$ with values in $Z(H)$. The identities in \eqref{bH} for $P \in A_1$ are 
the spectral differential equations in $x$ 
\[
L(x, \partial_x) \Psi_H(x,z) = \Psi_H(x,z) b(L)(z), \quad L(x, \partial_x) \in A_1.
\]
Let $A_2$ be the subalgebra of $B_2$ spanned by the differential-translation operators of the form $\partial_z^j T_z^{c H^{-1}}W$, $j \in \Nset$, 
$c \in \Cset$, $W \in Z(H)$. Let $K_1$ be the subalgebra of $B_1$ consisting of quasipolynomials in $z$ with values in $Z(H)$. 
The identities in \eqref{bH} for $P \in K_1$ are the spectral differential-translation equations in $z$ 
\[
b^{-1}(\La)(x)  \Psi_H(x,z) = \Psi_H(x,z) \cdot \La, \quad \La \in A_2.
\]

We are in the situation of \thref{GenTh} and can apply it to obtain functions that satisfy spectral differential equations in $x$ and 
spectral differential-translation equations in $z$ by Darboux transformations from $\Psi_H(x,z)$. 
This gives a second proof of the results of Kasman \cite{Ka}.

\bth{diff-tran} Let $H \in M_n(\Cset)$ be an invertible matrix.
If $L(x, \partial_x)$ is a $Z(H)$-valued differential operator with constant coefficients, 
and $Q(x, \partial_x)$ and $P(x, \partial_x)$ are $Z(H)$-valued differential operators
such that
\begin{equation}
\label{Hfact}
L(x, \partial_x) = Q(x, \partial_x) P(x, \partial_x),
\end{equation}
then 
\[
\Phi(x,z) := P(x, \partial_x) \Psi_H(x,z)
\]
is a generalized bispectral function, satisfying a spectral differential equation in $x$ and a spectral differential-translation 
equation in $z$.

More precisely, if 
\[
L(x, \partial_x) = Q'(x, \partial_x) g(x)^{-1} P'(x, \partial_x)
\]
for $Z(H)$-valued differential operators with quasipolynomial coefficients $Q'(x, \partial_x)$, $P'(x, \partial_x)$ 
and a nondegenerate $Z(H)$-valued quasipolynomial function $g(x)$,
then the function $\Phi(x,z) := P'(x, \partial_x) e^{xz} I_n$ satisfies
\begin{align*}
\left( P'(x, \partial_x) Q'(x, \partial_x) g(x)^{-1} \right) \Phi(x,z) &= \Phi(x,z) (b L)(z),
\\
g(x) \Phi(x,z) &= \Phi(x,z) \left( (bL)(z)^{-1} b(Q') b(P') \right).
\end{align*}
\eth
The first equation is a spectral differential equation in $x$. The second one is a spectral 
differential-translation equation in $z$; the point here is that $b(Q')$ and $b(P')$ are 
differential-translation operators. 

The classification of the factorizations of the form \eqref{Hfact} for monic differential operators $L$, $Q$ and $P$
is simpler than the classification in \thref{class-rank1}. The point is that the adelic type condition is not needed in this case
since the operators $P(x, \partial_x)$ and $Q(x, \partial_x)$ are not required to have rational coefficients. (For each 
such factorization $P(x, \partial_x)$ and $Q(x, \partial_x)$ have quasirational coefficients, recall \leref{lin} (a).)
Since $L(x,\partial_x)$ and $Q(x,\partial_x)$ are $Z(H)$-valued,
their vector kernels are nondegenerate $H$-invariant subspaces of $\QP_n$, recall \eqref{QPn}.
One easily shows that every nondegenerate $H$-invariant subspace of $\QP_n$
of dimension divisible by $n$ is the vector kernel of a unique monic differential operator $P(x, \partial_x)$ satisfying 
\eqref{Hfact} for some differential operators $L(x, \partial_x)$ and $Q(x, \partial_x)$ with the stated properties.
Thus we obtain the following:

{\em{The generalized bispectral Darboux transformations in \thref{diff-tran} for monic differential operators $L$, $Q$, and $P$ 
are classified by the nondegenerate $H$-invariant subspaces of $\QP_n$ of dimension divisible by $n$. For such a subspace, one constructs the 
unique monic differential operator $P(x, \partial_x)$ with that vector kernel and the corresponding generalized bispectral function $\Phi(x,z)$; 
this differential operator $P(x, \partial_x)$ is $Z(H)$-valued with quasipolynomial coefficients.}}
\sectionnew{Bispectral Darboux transformations from matrix Airy functions}
\label{Airy}
In this section we give an explicit classification of all bispectral Darboux 
transformations from the matrix Airy functions. We often 
make use of the spectral theory of matrix polynomials reviewed in the appendix. 
\subsection{}
\label{5.1}
We recall the setting of the (generalized) scalar Airy bispectral functions from \cite{BHY2}. 
Fix complex numbers $\alpha_0 \neq 0, \alpha_2, \ldots, \alpha_{N-1}$. Consider the scalar differential operator 
\[
M_{\Ai}(x, \partial_x) := \partial^N_x + \sum_{i=1}^{N-1}\alpha_i \partial^{N-i}_x + \alpha_0 x,
\]
called a (generalized) Airy operator. For $N=2$ this is the classical Airy operator up to normalization. 
It is well known that the operator $M_{\Ai}(x, \partial_x)$ has $N$-dimensional vector kernel over $\Cset$
consisting of entire functions.
Choose a basis of the kernel, $\psi_1(x), \ldots, \psi_N(x)$. 
For each function $\psi(x) \in \vker M$, define 
\begin{equation}
\label{xz-nota}
\psi(x,z) := \psi(x+z). 
\end{equation}
The latter satisfies the 
bispectral equations
\[
M_{\Ai}(x, \partial_x) \psi = z\psi \quad \text{and} \quad M_{\Ai}(z, \partial_z) \psi = x \psi.
\]
The bispectral Darboux transformations from it 
(in the sense of \thref{GenTh}) were classified in \cite{BHY2} and played an important role in the overall 
classification of scalar bispectral operators. 

Now we move to the matrix situation.
As in the matrix rank 1 case from Section \ref{rank1}, we take $\Om_1 = \Om_2 := \Cset$. Define the function
\[
\Psi_{\Ai}(x,z) = \psi(x,z) I_n
\] 
on $\Cset \times \Cset$, where $I_n$ is the identity matrix. Let $B_1$ and $B_2$ be both 
equal to the algebra of matrix-valued differential operators with polynomial coefficients. 
We will use the variables $x$ and $z$ for the first and second algebra, respectively.
Let $b \colon B_1 \to B_2$ be the isomorphism given by 
\[
b(x) := M_{\Ai}(z, \partial_z), \; \; b(\partial_x):= - \partial_z, \; \; b(W) := W, \; \forall \, W \in M_n(\Cset). 
\]
The first two equations imply that $b(M_{\Ai}(x, \partial_x)) = z$. We are in the bispectral situation of \prref{p1},
\[
P(x, \partial_x) \Psi_{\Ai}(x,z) = \Psi_{\Ai}(x,z) (bP)(z, \partial_z), \; \; \forall P(x, \partial_x) \in B_1. 
\]
Let $K_1$ and $K_2$ be the algebras of matrix-valued polynomial functions in $x$ and $z$, respectively. 
Denote by $A_1$ and $A_2$ the algebras consisting 
of operators of the form 
\[
q(M_{\Ai}(x, \partial_x)) \quad \mbox{and} \quad
q(M_{\Ai}(z, \partial_z))
\]
for matrix-valued polynomials $q(t)$. Then 
\[
b(A_1) = K_2, \; \; b(K_2) = A_1.
\]
The regular elements of the algebras $K_1$ and $K_2$ are the matrix-valued 
polynomial function $\theta(x)$ and $f(z)$ that are nondegenerate in the sense
that their determinants are nonzero polynomials.
By standard commutation arguments for matrix differential operators the sets 
$r(K_i)$ of regular elements are Ore subsets of $B_i$. 
The localizations $B_1[r(K_1)^{-1}]$ and $B_2[r(K_2)^{-1}]$
are isomorphic to the algebras of differential operators with rational coefficients 
in the variables $x$ and $z$. 
    
\thref{GenTh} implies the following:
\bth{Airy} Let $q(t)$ be a matrix-valued polynomial. Let 
$Q(x, \partial_x)$ and $P(x, \partial_x)$ be matrix-valued differential operators with rational coefficients, 
such that
\begin{equation}
\label{fact-Airy}
q(M_{\Ai}(x, \partial_x)) = Q(x, \partial_x) P(x, \partial_x).
\end{equation}
Then 
\[
\Phi(x,z) := P(x, \partial_x) \Psi_{\Ai}(x,z) 
\]
is a matrix bispectral function. More precisely, if 
\[
L(x, \partial_x) = Q'(x, \partial_x) g(x)^{-1} P'(x, \partial_x)
\]
for matrix-valued differential operators with polynomial coefficients $Q'(x, \partial_x)$, $P'(x, \partial_x)$ 
and a nondegenerate matrix-valued polynomial $g(x)$,
then $\Phi(x,z) := P'(x, \partial_x) \Psi_{\Ai}(x,z)$ satisfies
\begin{align*}
\left( P'(x, \partial_x) Q'(x, \partial_x) g(x)^{-1} \right) \Phi(x,z) &= \Phi(x,z) q(z),
\\
g(x) \Phi(x,z) &= \Phi(x,z) \left( q(z)^{-1} (bQ')(z, \partial_z) (bP')(z, \partial_z) \right).
\end{align*}
\eth
In \thref{Airy-class} below we classify all bispectral Darboux functions $\Phi(x,z)$ that enter 
in \thref{Airy} corresponding to $q(t)$ and $P(x, \partial_x)$ with nondegenerate 
leading terms.
\subsection{}
\label{5.2} 
First, we describe the kernels of the differential operators $L(x, \partial_x)$ in \thref{Airy}.
In the notation \eqref{xz-nota},
\[
\psi_i^{(j)} (x,\la) = (\partial_z^j \psi_i (x,z))|_{z=\la} = \partial_x^j \psi_i (x,\la), 
\quad \forall \, \la \in \Cset, 1 \leq i \leq N, j \in \Nset.
\]
Define the space
\[
\QA_n := \bigoplus_{\la \in \Cset, 1 \leq i \leq N, j \in \Nset} \, \psi_i^{(j)} (x,\la) \Cset^n.  
\]
Denote the standard basis of $\Cset^n$ by $\{ e_1, \ldots, e_n \}$. We have 
\[
\left( M_{\Ai}(x, \partial_x)  - \la \right) \psi_i^{(j)} (x,\la) = \psi_i^{(j-1)} (x,\la),
\]
and
\[
\{\psi_i^{(j)} (x,\la) \mid 1 \leq i \leq N, 0 \leq j \leq k-1 \}
\]
forms a basis of the kernel of $(M_{\Ai}(x, \partial_x) - \la)^k$ for all $\la \in \Cset$ 
and $k \in \Zset_+$, see \cite{BHY2}. One easily obtains from this that 
\begin{equation}
\{ \psi_i^{(j)} (x,\la) e_l \mid \la \in \Cset, 1 \leq i \leq N, j \in \Nset, 1 \leq l \leq n \}
\; \; 
\mbox{is a basis of} \; \; 
\QA_n.
\label{basis-QA}
\end{equation}
The fact that $\psi_i(x)$ are in the kernel of $M_{\Ai}(x, \partial_x)$ implies that 
\[
\QA_n = \bigoplus_{\la \in \Cset, 1 \leq i \leq N, 0 \leq j \leq N} \psi_i^{(j)} (x,\la) \Cset[x]^n
\]
and that
\begin{multline}
\{ \psi_i^{(j)} (x,\la) x^k e_l \mid \la \in \Cset, 1 \leq i \leq N, 0 \leq j \leq N-1, k \in \Nset, 1 \leq l \leq n \}
\\
\mbox{is a basis of} \; \; 
\QA_n.
\label{basis-QA2}
\end{multline}

Fix a monic matrix-valued polynomial   
\[
q(t) = \sum_{j=0}^d a_j t^j, \quad \mbox{where} \quad a_j \in M_n(\Cset), \; a_d = I_n.
\]
\bpr{ker-qM} For each matrix polynomial $q(t)$ as above, the vector kernel of the differential operator
$q(M_{\Ai}(x, \partial_x))$ over $\Cset$ is $ndN$-dimensional and consists of entire functions. 
More precisely, for each root $\la$ of $\det(q(t))$ and a choice of Jordan chains of $q(t)$ corresponding to $\la$ 
as in \thref{Jordan} $\{v_{0, l}, \ldots v_{k_l, l} \mid 1 \leq l \leq s \}$,  
the elements
\[
\sum_{r=0}^j \psi^{(r)}_i(x,\la) \frac{v_{j-r,l}}{r!}, \quad 1 \leq i \leq N, 0 \leq j \leq k_l
\] 
belong to $\vker q(M_{\Ai}(x, \partial_x))$. The set of all such elements for all roots 
$\la$ of $\det(q(t))$ and a corresponding set of Jordan chains as 
in \thref{Jordan} forms a basis of $\vker q(M_{\Ai}(x, \partial_x))$.
\epr
\begin{proof} Using the bispectrality of the functions $\psi_i(x,z)$ and the definition
of Jordan chains \eqref{App-J}, we obtain
\begin{align*}
&M_{\Ai}(x, \partial_x) \left( \sum_{r=0}^j \psi^{(r)}_i(x,\la) \frac{v_{j-r,l}}{r!} \right) = 
\sum_{r=0}^j \partial_z^j (M_{\Ai}(x, \partial_x)\psi_i (x,z))|_{z=\la}  \frac{v_{j-r,l}}{r!} \\
&= \sum_{r=0}^j \partial_z^j (\psi_i (x,z) q(z))|_{z=\la}  \frac{v_{j-r,l}}{r!} =
\sum_{r=0}^j \sum_{c=0}^r 
\begin{pmatrix}
r \\ c
\end{pmatrix} 
\psi^{(c)}_i (x,\la) q^{(r-c)}(\la)  \frac{v_{j-r,l}}{r!} = \\
&= \sum_{c=0}^j \left( \frac{\psi^{(c)}_i (x,\la)}{c!} 
\sum_{r=c}^j \frac{q^{(r-c)}(\la)}{(r-c)!} v_{(j-c)-(r-c),l} = 0 
\right).
\end{align*} 
To prove that the set in \prref{ker-qM} gives a basis of $\vker q(M_{\Ai}(x, \partial_x))$, 
first note that by \eqref{sum} the cardinality of this set equals $N \deg(\det(q(t)))$ = $ndN$. 
Because $\dim \vker q(M_{\Ai}(x, \partial_x)) \leq ndN$, 
all we need to show is that these elements are linearly independent. This follows from \eqref{basis-QA} 
and the linear independence in \thref{Jordan}. 
\end{proof} 
A direct argument with matrix polynomials leads to the following corollary of \prref{ker-qM}.  
\bco{LPQ-Airy} (a) Let $q(t)$ be a monic matrix polynomial.
If $Q(x, \partial_x)$ and $P(x, \partial_x)$ are monic matrix-valued differential operators 
with meromorphic coefficients in $\Cset$ which satisfy \eqref{fact-Airy},  
then $\vker P(x, \partial_x)$ is a nondegenrate subspace of $\QA_n$.

(b) For every nondenerate finite dimensional subspace $V$ of $\QA_n$ whose dimension is a multiple of $nN$, 
there exist a monic matrix polynomial $q(t)$ and monic matrix-valued differential 
operators $Q$ and $P$ satisfying the properties in part (a) 
such that $\vker P(x, \partial_x) = V$. 
\eco

Introduce a $\Zset_N$-action on $\QA_n$, where the generator $\sig$ of $\Zset_N$ acts 
by 
\[
\sig( \psi^{(j)}_i(x,\la) v):= \psi^{(j)}_{i+1}(x,\la) v, 
\quad \forall 
\; j \in \Nset, 1 \leq i \leq N, v \in \Cset^n
\]
and the $i+1$ term in the right hand side is taken modulo $N$ (i.e., $N+1 := 1$). Here we use 
the fact that the set in \eqref{basis-QA} is a basis of $\QA_n$.
In terms of the basis elements from \eqref{basis-QA2}, this action is given by
\[
\sig( \psi^{(j)}_i(x,\la)p(x)):= \psi^{(j)}_{i+1}(x,\la) p(x),
\quad \forall 
\; 0 \leq j \leq N-1, 1 \leq i \leq N, p(x) \in \Cset[x]^n.
\]
This follows at once from the fact that the functions $\psi_i(x)$ are in the kernel of $M_{\Ai}(x, \partial_x)$.
 
The next theorem is our classification result for all bispectral Darboux transformations 
from the matrix Airy functions. 
 
\bth{Airy-class} Assume that $q(t)$ is a monic matrix-valued polynomial,
and $Q(x, \partial_x)$ and $P(x, \partial_x)$ are monic matrix-valued differential operators whose coefficients are rational
functions. If \eqref{fact-Airy} is satisfied, then $\vker P(x, \partial_x)$ is a nondegenerate $\Zset_N$-invariant 
subspace of $\QA_n$ that has a basis consisting of elements of the form
\begin{equation}
\label{ker-P-Airy}
\sum_{j=0}^{N-1}\psi_i^{(j)} (x,\la) p_j(x),
\end{equation}
where $p_j(x) \in \Cset[x]^n$. The indices $i=1, \ldots, N$ and the complex numbers $\la \in \Cset$ 
can be different for different basis elements. 

Furthermore, for every such subspace of $\QA_n$ of dimension that is a multiple of $nN$, 
there exists a matrix polynomial $q(t)$ and a pair of operators 
$Q$ and $P$ with the above properties such that $\vker P(x, \partial_x)$ equals this subspace.
\eth
\bco{class-Airy} The bispectral Darboux transformations from the matrix Airy function $\Psi_{\Ai}(x,z)$ 
as in \thref{Airy} for $q(t)$ and $P(x, \partial_x)$ with nondegenerate leading terms are exactly
the functions of the form
\[
\Phi(x,z) = W(x) P(x, \partial_x) \Psi_{\Ai}(x,z),
\] 
where $P(x, \partial_x)$ is a monic matrix differential operator of order $dN$ with $ndN$-dimensional 
vector kernel which is a  
$\Zset_N$-invariant subspace of $\QA_n$ possessing a basis consisting of elements
of the form \eqref{ker-P-Airy} and $W(x)$ is a matrix-valued nondegenerate 
rational function. The operator $P(x, \partial_x)$ is uniquely reconstructed from the basis 
using the second part of \prref{fact}.
\eco
\subsection{}
\label{5.3}
We proceed with the proof of \thref{Airy-class}. The first part of the theorem 
follows from the following proposition and \coref{LPQ-Airy} (a).

\bpr{Airy-proof1} Assume that $P(x, \partial_x)$ is a matrix-valued differential operator 
with rational coefficients. If
\begin{equation}
\label{assum-Pker}
\sum_{\la \in \Cset} \sum_{i=1}^N \sum_{j=0}^{N-1} \psi_i^{(j)} (x,\la) p^{\la}_{j,i}(x) \in \vker P(x, \partial_x)
\end{equation}
for some vector-valued polynomials $p_{j,i}(x)$, then 
\[
\sig^k \left( \sum_{j=0}^{N-1} \psi_i^{(j)} (x,\la) p_{j,i}^\la(x) \right) \in \vker P(x, \partial_x), \quad
\forall \, \la \in \Cset, 1 \leq i \leq N, 0 \leq j, k \leq N-1.
\]
\epr
\begin{proof} It follows from \eqref{basis-QA2} that 
\[
\{ \psi_i^{(j)} (x,\la) e_l \mid \la \in \Cset, 1 \leq i \leq N, 0 \leq j \leq N-1,  
1 \leq l \leq n \}
\]
is a linearly independent set over $\Cset(x)$. Thus, \eqref{assum-Pker} and the assumption that 
$P(x, \partial_x)$ has rational coefficients imply that 
\begin{equation}
\label{inst}
\sum_{j=0}^{N-1} \psi_i^{(j)} (x,\la) p^{\la}_{j,i}(x) \in \vker P(x, \partial_x), \quad
\forall \, \la \in \Cset, 1 \leq i \leq N.
\end{equation}
Using the facts that each such element belongs to $\vker P(x, \partial_x)$ and 
$P(x, \partial_x)$ has rational coefficients, and recursively expressing 
\[
\psi_i^{(N)} (x,\la) = -\sum_{j=0}^{N-1}\alpha_i \psi_i^{(j)}(x, \la) - \alpha_0 (x+ \la) \psi_i (x, \la),
\]
produces a $\Cset[x]$-linear combination of 
\[
\{ \psi_i^{(j)} (x,\la) e_l \mid 0 \leq j \leq N-1, 1 \leq l \leq n \}
\]
that equals 0. In view of \eqref{basis-QA2}, this is only possible if each coefficient of the combination equals 0. 
However, exactly the same thing would happen if we apply $P(x, \partial_x)$ to
\[
\sum_{j=0}^{N-1} \psi_{i+k}^{(j)} (x,\la) p^{\la}_{j,i}(x) = 
\sig^k \left( \sum_{j=0}^{N-1} \psi_i^{(j)} (x,\la) p^{\la}_{j,i}(x) \right)
\]
instead of \eqref{inst}. Thus, for all $0 \leq k \leq N-1$, 
these elements belong to $\vker P(x, \partial_x)$ 
which completes the proof of the proposition.
\end{proof}

The second part of \thref{Airy-class} follows from the following proposition and \coref{LPQ-Airy} (b).
\bpr{Airy-proof2} Assume that $V$ is a nondegenerate $ndN$-dimensional $\Zset_N$-invariant subspace of $\QA_n$ having a basis 
consisting of elements of the form \eqref{ker-P-Airy}, where $p_j(x) \in \Cset[x]^n$. (The indices $i=1, \ldots, N$ and the complex 
numbers $\la \in \Cset$ can be different for different basis elements.) Let $P(x, \partial_x)$ be 
the unique monic matrix-valued differential operator on $\Cset$ of order $dN$ whose vector kernel kernel 
equals $V$ (recall \prref{diff-op} (b)). Then $P(x, \partial_x)$ has rational coefficients.  
\epr
\begin{proof} Denote the operator
\[
P(x, \partial_x) = \partial_x^{dN} + \sum_{l=0}^{dN-1} a_l(x) \partial_x^l,
\]
where $a_l(x)$ are matrix-valued meromorphic functions on $\Cset$.

The definition of the $\Zset_N$-action on $\QA_n$ implies that we can assume 
that the basis of $V$, consisting of elements of the form \eqref{ker-P-Airy}, is itself $\Zset_N$-invariant.
Let 
\begin{equation}
\label{Nelem}
\sig^k \left( \sum_{j=0}^{N-1}\psi_i^{(j)} (x,\la) p_j(x) \right)
= \sum_{j=0}^{N-1}\psi_{i+k}^{(j)} (x,\la) p_j(x), \quad 1 \leq k \leq N-1
\end{equation}
be $N$ elements in the vector kernel of $P(x, \partial_x)$. Using the condition that 
they belong to $\vker P(x, \partial_x)$ and expressing 
\[
\psi_i^{(N)} (x,\la) = -\sum_{j=1}^{N-1}\alpha_i \psi_i^{(j)}(x, \la) - \alpha_0 (x+ \la) \psi_i (x, \la),
\]
produces a $\Cset[x]$-linear combination of 
\[
\{ \psi_i^{(j)} (x,\la) e_l \mid 0 \leq j \leq N-1, 1 \leq l \leq n \}
\]
that equals 0. Once again \eqref{basis-QA2} implies that each coefficient 
of the combination should equal 0. Thus, the fact that one of the elements 
\eqref{Nelem} belongs to $\vker P(x, \partial_x)$ 
is equivalent to imposing $N$ conditions of the form
\[
\sum_{l=0}^{dN-1} a_l(x) c_{l, j}(x) = 0, \quad \forall \, 0 \leq j \leq N-1,  
\]
where $c_{l, j}(x)$ are some vector-valued polynomial functions uniquely determined from 
$p_0(x), \ldots, p_{N-1}(x)$. However, the condition that every element in \eqref{Nelem} belongs to $\vker P(x, \partial_x)$ 
is equivalent to imposing exactly the same $N$ conditions because $\psi_i(x, \la)$ satisfy the
same Airy equation. 

So, the conditions on the operator $P(x, \partial_x)$ (having vector kernel equal to $V$) 
are equivalent to imposing $dN$ linear vector conditions 
on the matrix meromorphic functions $\{a_l(x) \mid 0 \leq l \leq dN-1 \}$. Because the coefficients of these 
conditions are vector-valued polynomial functions and there is a unique monic differential operator $P(x, \partial_x)$ 
with the stated properties, all its coefficients $a_l(x)$ should be rational.  
\end{proof}
\sectionnew{Examples}
\label{exe}
In this section we present several examples illustrating the classification results from the previous two sections.
\subsection{}
\label{6.1}
First, we give three examples of matrix rank 1 bispectral Darboux transformations, classified in \thref{class-rank1}
and \coref{class-rank1-b}.

\bex{ex-1} In the case of $2\times 2$ matrices, consider the differential operator
\[
L(x, \partial_x) := \partial_x^2 I_2 \in A_1.
\]
(Below all instances of the identity matrix $I_n$ which are clear from the setting will be suppressed.) 
The vector kernel of $L$ consists of all vector linear functions. In the setting of \thref{class-rank1}, 
take the subspace $V$ of $\QP_2$ with basis  
\[
f_1(x) := (x , 0)^t, \quad f_2(x) = (a, x)^t
\]
for an arbitrary $a \in \Cset$. Set
\[
F_1(x) := (f_1(x), f_2(x)) = 
\begin{pmatrix}
x & a \\
0 & x
\end{pmatrix}.
\]
Consider the corresponding operator $P(x, \partial_x)$ given by \eqref{Pf}, which in this case simplifies 
to 
\[
P(x, \partial_x) = \partial_x - F_1'(x)F_1^{-1}(x) = \partial_x - \begin{pmatrix}
x^{-1} & -ax^{-2}  \\ 
0 &x^{-1}\\ 
\end{pmatrix}.
\]
By \thref{class-rank1}, we have the factorization $L(x, \partial_x) = Q(x, \partial_x) P(x, \partial_x)$ 
from which the operator $Q$ is computed to be
\[
Q(x, \partial_x) = \partial_x + \begin{pmatrix}
x^{-1} & - ax^{-2}  \\ 
0 &x^{-1}\\ 
\end{pmatrix}.
\]
This gives rise to the bispectral Darboux transformation
\[
\Phi(x, z) := P(x, \partial_x) (e^{xz} I_2) = 
e^{xz} \begin{pmatrix}
z- x^{-1} & a x^{-2}  \\ 
0 & z- x^{-1}
\end{pmatrix}.
\]
It satisfies 
\[
\wt{L}(x, \partial_x) \Phi(x, z) = \Phi(x, z) z^2
\]
where
\[
\wt{L}(x, \partial_x) := P(x, \partial_x) Q(x, \partial_x)
  = \partial^2  -  \begin{pmatrix}
2x^{-2} &  - 4ax^{-3}  \\ 
0 &2x^{-2}\\ 
\end{pmatrix}.    
\]
For the dual spectral equation, given by \thref{rank1}, we need to represent $L(x, \partial_x)$ 
in the form \eqref{special-fact}:
\[
\partial_x^2 I_2 = Q'(x, \partial_x) \frac{1}{x^4} P'(x, \partial_x) \; \; \mbox{with} \; \; 
Q' := \partial_x x^2+ \begin{pmatrix}
x & - a   \\ 
0 &x\\ 
\end{pmatrix}, \; 
P' := x^2\partial_x  -   \begin{pmatrix}
x & - a \\ 
0 &x
\end{pmatrix} \in B_1.
\]
Now, by \thref{rank1}, the function $\Phi(x,z)$ satisfies the dual spectral equation
\[
\Phi(x,z) \wt{\Lambda}(z, \partial_z) = x^4 \Phi(x,z)
\]
where
\begin{align*}
\wt{\Lambda}(z, \partial_z) &:= 
\frac{1}{z^2} b \left[ \partial_x x^2+ \begin{pmatrix}
x & - a   \\ 
0 &x\\ 
\end{pmatrix} 
\right] 
\circ
b \left[ x^2 \partial_x  -  \begin{pmatrix}
x & - a \\ 
0 &x
\end{pmatrix}       
\right]
\\
&= 
\frac{1}{z^2} \left[z \partial_z^2 - \partial_z + \begin{pmatrix}
0 & - a   \\ 
0 &0\\ 
\end{pmatrix} \right] 
\circ
\left[ \partial_z^2 z  + \partial_z + \begin{pmatrix}
0 & a \\ 
0 & 0
\end{pmatrix}       
\right]
\\
&= z^{-2} ( z \partial_z^2 - \partial_z)( \partial_z^2 z + \partial_z) I_2  + 
\begin{pmatrix}
0 & - 2 a  z^{-2} \partial_z  \\ 
0 &0\\ 
\end{pmatrix}.
\end{align*}
When $a=1$, this is the first example in \cite{BL} and \cite{G1}, where 
a (slightly more complicated) third order dual spectral equation was considered.  
\eex

\bex{ex-2} For $3\times 3$ matrices, consider the operator $L = \partial_x^2 I_3 \in A_1$. Its vector kernel consists 
of vector linear functions. Let $V$ be the subspace of $\QP_3$ with a basis consisting of the 
columns of the matrix  
\[
F_1(x) :=
\begin{pmatrix}
x&  1 & 0 \\ 
0 &x &1 \\ 
0&0  &x
\end{pmatrix}.
\]
The corresponding operator $P(x, \partial_x)$, given by \eqref{Pf}, is
\[
P(x, \partial_x) := \partial_x - F_1'(x) F_1(x)^{-1} = 
\partial_x - (x^{-1} I_3 - x^{-2} J + x^{-3} J^2)
\]
where $J$ is the super-diagonal matrix
\[
J:=
\begin{pmatrix}
0 & 1 & 0 \\ 
0 & 0 & 1 \\ 
0 & 0 & 0
\end{pmatrix}.
\]
By \thref{class-rank1}, we have the factorization $L(x, \partial_x) = Q(x, \partial_x) P(x, \partial_x)$ 
from which the operator $Q$ is computed to be
\[
Q(x, \partial_x) = \partial_x + (x^{-1} I_3 - x^{-2} J + x^{-3} J^2). 
\]
The bispectral Darboux transformation
\[
\Phi(x, z) := P(x, \partial_x) (e^{xz} I_3) = 
e^{xz} \begin{pmatrix}
z - x^{-1} &  x^{-2}     & -x^{-3}\\ 
0          & z - x^{-1}  &  x^{-2} \\ 
0          & 0           &  z - x^{-1}
\end{pmatrix}
\]
satisfies 
\[
\wt{L}(x, \partial_x) \Phi(x, z) = \Phi(x, z) z^2
\]
where
\[
\wt{L}(x, \partial_x) := P(x, \partial_x) Q(x, \partial_x) =
\partial_x^2 - 2 x^{-2} I_3 + 4 x^{-3} J - 6 x^{-4} J^2.
\]
To write the dual spectral equation in \thref{rank1}, we represent $L(x, \partial_x)$ 
in the form \eqref{special-fact}, $\partial_x^2 I_2 = Q'(x, \partial_x) x^{-3} P'(x, \partial_x)$
with
\[
Q'(x, \partial_x) := \partial_x x^3 + (x^2 I_3 - x J + J^2),
\quad
P'(x, \partial_x) := x^3 \partial_x - (x^2 I_3 - x J + J^2)
\in B_1.
\]
\thref{rank1} implies that $\Phi(x,z)$ satisfies the dual spectral equation
\[
\Phi(x,z) \wt{\Lambda}(z, \partial_z) = x^4 \Phi(x,z)
\]
with
\begin{align*}
\wt{\Lambda}(z, \partial_z) &:= 
z^{-2} b(Q')(z, \partial_z) b(P')(z, \partial_z) 
\\
&= z^{-2} (z \partial_z^3 I_3 - (\partial_z^2 I_3 + \partial_z J + J^2))
(\partial_z^3 z I_3 + (\partial_z^2 I_3 + \partial_z J + J^2)).
\end{align*}
This bispectral function is the second example in \cite{G1}, where a conjecture is stated about the 
full algebra of dual differential operators with eigenfunction $\Phi(x,z)$. It contains lower order 
operators, e.g., the degenerate operator 
\[
\left(\begin{matrix}
0 &     0  & 0 \\ 
0 & 0  &  0 \\ 
1 & 0 &0
\end{matrix}       
\right) \partial_z^2 +   \left(\begin{matrix}
0 &     0  & 0 \\ 
1 & 0  &  0 \\ 
-2 z^{-1}& 1&0
\end{matrix}       
\right) \partial_z  +   \left(\begin{matrix}
1 &     0  & 0 \\ 
-2 z^{-1} & 2  &  0 \\ 
0 & 0 &1
\end{matrix}       
\right).  
\]
\eex

\bex{ex-3} For $3\times 3$ matrices, consider the operator
\[
L(x, \partial_x) = \partial_x^2 I_3 + \partial_x J  + J^2 = 
\begin{pmatrix} 
\partial_x^2& \partial_x  & 1\\
0 & \partial_x^2 & \partial_x \\
0 & 0& \partial_x^2
\end{pmatrix} \in A_1.
\]
A basis of $\vker L(x, \partial_x)$ is given by
\[
(1, 0, 0)^t, (x, 0, 0)^t,
(0, 1, 0)^t,  (-x^2/2, x, 0)^t, (0,-x^2/2,x)^t,
(-x^2/2, 0, 1)^t.
\]
Consider the subspace of $\QP_3$ with a basis consisting of the columns of the matrix function
\[
F_1(x) =
\begin{pmatrix} 
x & -x^2/2 & -x^2/2\\
0 & x   & 0\\
0 & 0& 1
\end{pmatrix}.
\]
The corresponding operator $P(x, \partial_x)$, given by \eqref{Pf}, is  $P(x, \partial_x) := \partial_x - F_1'(x) F_1(x)^{-1}$, or explicitly,
\[
P(x, \partial_x)=
\partial_x -
\begin{pmatrix} 
x^{-1} & -3/2 & -x/2\\
0 & x^{-1} & 0 \\
0 & 0& 0
\end{pmatrix} =x^{-1} \Big\{ x \partial_x -
\begin{pmatrix} 
1 & -3x/2 & -x^2/2\\
0 & 1 & 0 \\
0 & 0& 0
\end{pmatrix}    \Big \} = x^{-1} P'.
\]
By \thref{class-rank1}, we have the factorization $L(x, \partial_x) = Q(x, \partial_x) P(x, \partial_x)$ 
from which the operator $Q$ is computed to be
 \[
 Q = \partial_x + \begin{pmatrix} 
 x^{-1} & -1/2   & -x/2\\
 0 & x^{-1} & 1 \\
 0 & 0& 0
 \end{pmatrix} = \Big\{   \partial_x x + \begin{pmatrix} 
 1 & -x/2 & -x^2/2\\
 0 & 1 & x \\
 0 & 0& 0
 \end{pmatrix}  \Big\}x^{-1} =  Q' x^{-1} .
 \]
 The operators $P'(x,\partial_x)$ and $Q'(x,\partial_x)$ will be needed for the dual bispectral equation.
 
 The bispectral Darboux transformation
\[
\Phi(x,z) := P(x, \partial_x) (e^{xz} I_3) = 
\begin{pmatrix} 
z- x^{-1} & 3/2 & x/2\\
0 & z- x^{-1} & 0 \\
0 & 0& z
\end{pmatrix}
\]
satisfies 
\[
\wt{L} (x,\partial_x)\Phi(x,z) = \Phi(x,z) \begin{pmatrix} 
z^2& z  & 1\\
0 & z^2 & z \\
0 & 0& z^2 
\end{pmatrix} 
\]
where the differential operator $\wt{L}(x,\partial_x)$ is given by
\[
\wt{L}(x,\partial_x):= P(x,\partial_x) Q(x,\partial_x),
\]
omitting the detailed expansion for the sake of brevity. By \thref{rank1}, the function $\Phi(x,z)$ satisfies the 
dual spectral equation
\[
x^2 \Phi(x,z)  = \Phi(x,z) \wt{\La}(x, \partial_x)
\]
where 
\begin{align*}
&\wt{\Lambda}(z, \partial_z) = (z^2 I_3 + z J + J^2)^{-1} b(Q')(z, \partial_x) b(P') (z, \partial_z) =  \\
&
\begin{pmatrix}
z^{-2} & -z^{-3} & 0 \\
0 &  z^{-2} & - z^{-3} \\
0 & 0 & z^{-2}
\end{pmatrix}
\Big\{  z \partial_z - \begin{pmatrix} 
 1 & \partial_z/2 & -\partial_z^2/2\\
 0 & 1 & - \partial_z \\
 0 & 0& 0
 \end{pmatrix}  \Big\}
\Big\{ \partial_z z +
\begin{pmatrix} 
1 & 3 \partial_z/2 & -\partial_z^2/2\\
0 & 1 & 0 \\
0 & 0& 0
\end{pmatrix}  
\Big\}.
\end{align*}
\eex

\subsection{}
\label{6.2}
We finish with an example illustrating the classification results for all bispectral Darboux transformations 
from the matrix Airy functions. 

\bex{A-1} Consider the matrix polynomial
\[
q(t) = (I_2 t - J)^2 = I_2 t^2 - 2 J t 
\]
where 
\[
J := \begin{pmatrix} 
0 & 1\\
0 & 0
\end{pmatrix}.
\]
For the second order Airy operator $M_{\Ai}(x,\partial_x) = \partial_x^2 - x$,
consider the operator 
\[
L(x,\partial_x) := q (M_{\Ai}(x,\partial_x)) = (I_2 \partial_x^2 - I_2 x - J)^2 \in B_1.
\]
It satisfies $b(L(x, \partial_x))= q(z) = z^2 I_2 - 2z J$, i.e.,  
\[
L(x, \partial_x) \Psi_{\Ai}(x,z) = \Psi_{\Ai}(x,z) (z^2 I_2 - 2z J). 
\]

According to our scheme (from \prref{ker-qM}, \thref{Airy-class} and \coref{class-Airy}), we first need to find the Jordan chains 
of the monic matrix polynomial $q(t)$.
Its characteristic polynomial $\chi(t) = \det(q(t)) = t^4$ has only one root 0 and
\[
q(0) = 0,  \quad 
q'(0) = - 2 J, \quad q''(0) = 2 I_2.
\]
Denote by $e_1, e_2$ the standard basis of $\Cset^2$. Then
\begin{align*}
& e_2, 
\\
& e_1, (1/2) e_2, e_1
\end{align*}
is a maximal family of Jordan chains as in \thref{Jordan}. 
Recall from \S \ref{5.1} that $\psi_1(x), \psi_2(x)$ denotes a basis of $\vker M_{\Ai}(x,\partial_x)$.
\prref{ker-qM} produces the following basis of $\vker L(x, \partial_x):$
\[
\begin{pmatrix}
0 \\
\psi_i(x)
\end{pmatrix}, \; 
\begin{pmatrix}
\psi_i(x) \\
0
\end{pmatrix}, \;
\begin{pmatrix}
\psi'_i(x) \\
\psi_i(x)/2
\end{pmatrix}, \;
\begin{pmatrix}
\psi_i(x) + \psi''_i(x)/2 \\
\psi'_i(x)/2
\end{pmatrix},
\quad i=1,2.
\]
(One can further simplify it by taking linear combinations.)
Consider the 4-dimensional subspace $V$ of $\QA_2$ with basis 
\[
\begin{pmatrix}
\psi'_i(x) + \psi_i(x) \\
\psi_i(x) 
\end{pmatrix}, \; 
\begin{pmatrix}
\psi''_i(x) \\
\psi'_i(x)
\end{pmatrix}, \quad
i = 1,2.
\]
It satisfies the conditions in \thref{Airy-class}, in particular $V$ is $\Zset_2$-invariant.   It can be computed by using quasideterminants.
However, more effectively, it can be computed by setting 
$P(x, \partial_x) = \partial_x^2 + b(x) \partial_x + c(x)$ 
and then solving for $b(x)$ and $c(x)$ from the equations $P(x,\partial_x) f(x) = 0$ for the 4 basis elements $f(x)$ of $V$.
The concrete solution of the equations for $b(x)$ and $c(x)$ is done by expressing all derivatives $\psi^{(k)}_i(x)$ in terms of 
$\psi_i(x)$ and $\psi'_i(x)$ using the Airy equation, and then setting the coefficients in front of $\psi_i(x)$ and $\psi'_i(x)$
in the equations $P(x,\partial_x) f(x) = 0$ to be equal to 0. This gives that
\[ 
	P(x, \partial_x) = \partial_x^2 +  
	\frac{1}{x-1}
	\begin{pmatrix}
		0 & 1- x  
		\\
		1 & -1 
	\end{pmatrix}
	\partial_x - 
	\frac{1}{x-1}
	\begin{pmatrix}
		  (x- 1)^2  &   2x - 2
		\\
		0 & x^2
	\end{pmatrix}.
\]
By \thref{class-rank1}, $L(x, \partial_x) = Q(x, \partial_x) P(x, \partial_x)$ for some second order monic 
differential operator $Q(x,\partial_x)$. By comparing coefficients one obtains
\[
	Q(x, \partial_x) = \partial^2_x  -  \partial_x \frac{1}{x-1} 
	\begin{pmatrix}
		0 & 1 - x \\
			 1 & - 1
	\end{pmatrix}
	  + \frac{1}{(x-1)} 
	\begin{pmatrix}
	- x^2  & 1 \\
          0      & - (x-1)^2
	\end{pmatrix}.
\]
The bispectral Darboux transformation
\[
\Phi(x, z) := P(x, \partial_x) \Psi_{\Ai}(x,z)
\]
satisfies 
\[
\wt{L}(x, \partial_x) \Phi(x, z) = \Phi(x, z) (z^2 I_2 - 2 J)
\]
where
\[
\wt{L}(x, \partial_x) := P(x, \partial_x) Q(x, \partial_x).
\]
For the sake of brevity we leave the detailed computation of $\wt{L}(x, \partial_x)$ to the reader.

To write the dual spectral equation from \thref{Airy}, we represent $q(M_{\Ai}(x,\partial_x))$ 
in the form \eqref{special-fact}, 
\[
q(M_{\Ai}(x,\partial_x)) = Q'(x, \partial_x) \frac{1}{(x- 1)^2} P'(x, \partial_x)
\]
with
\[
Q'(x, \partial_x) := Q(x,\partial_x)(x-1),
\quad
P'(x, \partial_x) := (x-1) P(x, \partial_x) 
\in B_1.
\]
Explicitly we have
\begin{align*}
		P'(x, \partial_x) &=  (x-1)\partial_x^2 +  
		\begin{pmatrix}
		0 & 1- x  
		\\
		1 & -1 
		\end{pmatrix}
		\partial_x - 
		\begin{pmatrix}
		(x- 1)^2  &   2x - 2
		\\
		0 & x^2
		\end{pmatrix},
		\\
		Q'(x, \partial_x) & = \partial^2_x (x-1)-  \partial_x 
	\begin{pmatrix}
		0 & 1 - x \\
			 1 & - 1
	\end{pmatrix}
	  +
	\begin{pmatrix}
	- x^2  & 1 \\
          0      & - (x-1)^2
	\end{pmatrix}.
\end{align*}
\thref{Airy} implies that $\Phi(x,z)$ satisfies the dual spectral equation
\[
\Phi(x,z) \wt{\Lambda}(z, \partial_z) =(x-1)^3 \Phi(x,z)
\]
with
\begin{multline*}
\wt{\Lambda}(z, \partial_z) := (z I_2 -J)^{-2} b(Q')(z, \partial_z) b(P')(z, \partial_z) \\
= (z^{-2} + 2 z^{-1} J) b(Q')(z, \partial_z) b(P')(z, \partial_z) 
\end{multline*}
where
\begin{eqnarray*}
		b(P')(z, \partial_z)& = & (M_z -1)\partial_z^2 -  
		\begin{pmatrix}
		0 & 1- M_z 
		\\
		1 & -1 
		\end{pmatrix}
		\partial_z - 
		\begin{pmatrix}
		(M_z- 1)^2  &   2 M_z - 2
		\\
		0 & M_z^2
		\end{pmatrix},
		\\
		b(Q')(z, \partial_z) &=& \partial_z^2 (M_z-1)^2 +  
		\partial_z \begin{pmatrix}
		0 & 1 - M_z \\
		1 & - 1
		\end{pmatrix}
		 +  \begin{pmatrix}
	             - M_z^2  & 1 \\
                     0      & - (M_z-1)^2
	         \end{pmatrix}
\end{eqnarray*}
and we abbreviated $M_x :=M_{\Ai}(x, \partial_x)$, $M_z :=M_{\Ai}(z, \partial_z)$.
\eex
\sectionnew{Appendix. Matrix polynomials}
\label{Appen}

Here we collect some facts about the spectral theory of matrix polynomials, used in the paper. We follow \cite{GLR}
and refer the reader to \cite[Sections 1.4-1.6 and S1.5]{GLR} for details.

Consider a monic matrix polynomial
\[
q(t) = \sum_{j=0}^d a_j t^j, \quad \mbox{where} \quad a_j \in M_n(\Cset), \; a_d = I_n.
\]
Define the monic polynomial 
\[
\chi(t) := \det (q(t))
\]
of order $nd$. It equals the characteristic polynomial of the companion matrix
\begin{equation}
C = \left( \begin{matrix}
0 & I & 0 & \ldots & 0 \\
0 & 0 & I &\ldots & 0 \\
\vdots & \vdots & \vdots & \vdots & \vdots \\
- a_0 & - a_1 & -a_2& \ldots & - a_{d-1}
\end{matrix} 
\right) \in M_{nd}(\Cset),
\label{comp-M}
\end{equation}
see \cite[Theorem 1.1]{GLR}.

A sequence of vectors $v_0, v_1  \ldots , v_k \in \Cset^n$, with $v_0 \neq 0$ is called a {\em{Jordan chain}} of $q(t)$ 
of length $k+1$ corresponding to $\lambda \in \Cset$ if the following equalities holds:
\begin{equation}
\sum_{r=0}^{j} \frac{1}{r!} q^{(r)}(\lambda)v_{j-r} = 0,\; \forall \, 0 \leq j \leq  k. 
\label{App-J}
\end{equation}
Such a chain exists only if $\la$ is a root of $\chi(t)$. Its leading term $v_0$ is an eigenvector of $q(\lambda)$.
The condition \eqref{App-J} is equivalent to saying that the function
\[
u(x):= \left( \frac{x^k}{k!} v_0 + \frac{x^{k-1}}{(k-1)!} v_1 + \cdots + v_k \right) e^{\lambda x} 
\]
is a solution of the ODE $q(\partial_x) u(x) =0$. When $q(t)$ is the characteristic polynomial 
of a square matrix $A$, this notion recovers the classical notion of a Jordan chain for $A$.

The following theorem is a generalization of the well known Jordan normal form theorem for 
square matrices, see \cite[Theorem 1.12]{GLR} for details.

\bth{Jordan} Let $q(t)$ be a monic matrix polynomial and $\la$ be a root of $\chi(t)$ of multiplicity 
$m$. Then there exist positive integers $k_1, \ldots, k_s$ and Jordan chains
\[
v_{0,l}, \ldots, v_{k_l, l} \in \Cset^n
\] 
(for $1 \leq l \leq s$) of $q(t)$ with respect to $\la$ such that $\{v_{0, 1}, \ldots, v_{0,s} \}$ 
is a basis of $\ker q(\la)$ and 
\begin{equation}
\label{sum}
(k_1+1) + \cdots + (k_s +1) = m.
\end{equation}
\eth
The lengths of the Jordan chains in the theorem can be read off from the local Smith form of $q(t)$
in terms of the partial multiplicities of $q(t)$ at $\la$, see 
\cite[Section S1.5]{GLR}. 

\end{document}